\documentclass[final,oneside]{amsart}
\usepackage{geometry}
\usepackage[utf8]{inputenc}
\usepackage[verbose]{hyperref}
\usepackage[english]{babel}
\usepackage{enumerate}
\usepackage{aliascnt}
\usepackage{amssymb}
\usepackage{amsmath}
\usepackage{verbatim}
\usepackage{mathtools}
\usepackage{tikz}
\usepackage[notcite,notref]{showkeys}
\usepackage{amsthm}
\usepackage[title]{appendix}

\numberwithin{equation}{section}

\usepackage{color}

\usepackage{dsfont} 
\usepackage{hyperref} 
\usepackage{cleveref}

\allowdisplaybreaks[4] 

\newtheorem{prop}{Proposition}[section]

\newaliascnt{lem}{prop} 

\newtheorem{lem}[lem]{Lemma}

\aliascntresetthe{lem} 
\Crefname{lem}{Lemma}{Lemmas}

\newaliascnt{defi}{prop} 
 \newtheorem{defi}[defi]{Definition}

\aliascntresetthe{defi} 
\Crefname{defi}{Definition}{Definitions}

\newaliascnt{cor}{prop} 
 \newtheorem{cor}[cor]{Corollary}
 
\aliascntresetthe{cor}

\newaliascnt{remark}{prop} 
 \newtheorem{remark}[remark]{Remark}
 
\aliascntresetthe{remark} 

\newaliascnt{thm}{prop} 
 \newtheorem{thm}[thm]{Theorem}
 
\aliascntresetthe{thm} 

\newaliascnt{example}{prop} 
 \newtheorem{example}[example]{Example}
 
\aliascntresetthe{example}

\def\equationautorefname~#1\null{%
  (#1)\null
}

\newcommand\mydots{\hbox to 1em{.\hss.\hss.}} 

\newcommand{\blue}[1]{#1}

\newcommand{\R}{\ensuremath{\mathbb{R}}}
\newcommand{\N}{\ensuremath{\mathbb{N}}}

\newcommand*\diff{\mathop{}\!\mathrm{d}}

\newcommand{\defeq}{\vcentcolon=}
\newcommand{\eqdef}{=\vcentcolon}
\newcommand{\CalE}{\ensuremath{\mathcal{E}}}
\newcommand{\CalC}{\ensuremath{\mathcal{C}}}
\newcommand{\CalD}{\ensuremath{\mathcal{D}}}
\newcommand{\CalL}{\ensuremath{\mathcal{L}}}
\newcommand{\CalA}{\ensuremath{\mathcal{A}}}
\newcommand{\CalF}{\ensuremath{\mathcal{F}}}
\newcommand{\X}{\ensuremath{\mathbb{X}}}
\newcommand{\Y}{\ensuremath{\mathbb{Y}}}
\newcommand{\vKap}{\ensuremath{\vec{\kappa}}}
\newcommand{\abs}[1]{\ensuremath{\lvert #1\rvert}}

\newcommand{\Norm}[2]{\ensuremath{\left\Vert #1 \right\Vert_{#2}}}
\newcommand{\norm}[2]{\ensuremath{\Vert #1 \Vert_{#2}}}

\newcommand{\dtzero}{\left.\frac{\diff}{\diff t}\right\vert_{t=0}}

\DeclareMathOperator{\Id}{Id}

\DeclareMathOperator{\tr}{tr}

\title[Existence and convergence for the length-preserving elastic flow]{Existence and convergence of the length-preserving elastic flow of clamped curves}
\author[F.~Rupp]{Fabian Rupp}
\address[F.~Rupp]{Faculty of Mathematics, University of Vienna, Oskar-Morgenstern-Platz 1, 1090 Vienna, Austria.} \email{fabian.rupp@univie.ac.at}
\author[A.~Spener]{Adrian Spener}
\address[A.~Spener]{Institute of Applied Analysis, Ulm University, Helmholtzstra\ss e 18, 89081 Ulm, Germany.}
\email{adrian.spener@alumni.uni-ulm.de}

\begin{document}

\begin{abstract}
We study the evolution of curves with fixed length and clamped boundary conditions moving by the negative $L^2$-gradient flow of the elastic energy.
For any initial curve lying merely in the energy space we show existence and parabolic smoothing of the solution. Applying previous results on long time existence and proving a constrained {\L}ojasiewicz--Simon gradient inequality we furthermore show convergence to a critical point as time tends to infinity.
\end{abstract}

\maketitle

 
 
 

\section{Introduction and main results}

For an immersed curve $f\colon I\defeq [0,1]\to \R^d$, $d\geq 2$, its \emph{Euler--Bernoulli energy} or simply \emph{elastic energy} is defined by
\begin{align*}
	\CalE(f)\defeq \frac{1}{2}\int_{I} \abs{\vKap}^2\diff s.
\end{align*}
Here $\diff s\defeq \gamma\diff x$, where $\gamma\defeq \abs{\partial_x f}$ denotes the \emph{arc-length element}, and $\vKap \defeq \partial_s^2 f$ is the \emph{curvature vector field}, where $\partial_s \defeq \gamma^{-1}\partial_x$ is the \emph{arc-length derivative}. 

In this article, we deform an initial curve $f_0$ in such a way that its elastic energy decreases as fast as possible, while keeping the \emph{(total) length} $\CalL(f)\defeq \int_I\diff s$ fixed. This yields the geometric evolution equation
\begin{align}\label{eq:EFEvolution}
	\partial_t^{\perp}f = - \nabla_s^2\vKap-\frac{1}{2}\abs{\vKap}^2\vKap+\lambda \vKap.
\end{align}
Here $\nabla_s$ denotes the connection on the normal bundle along $f$, i.e.\ $\nabla_s \defeq P^{\perp}\partial_s$, where 
$P^{\perp}X \defeq X^{\perp_f}\defeq X-\langle X,\partial_s f\rangle\partial_s f$ denotes the orthogonal projection along $f$ of any {vector field $X$ along $f$}. The  \emph{Lagrange multiplier}  $\lambda$ depends on the solution $f$ and is given by 
\begin{equation}
\label{eq:lambda}
\lambda(f) = \lambda(f)(t)= \frac{\int_I \left\langle 
	\nabla_s^2 \vKap + \frac{1}{2} \abs{\vKap} ^2 \vKap, \vKap \right\rangle \diff s}{\int_I \vert \vKap \vert^{^2} \diff s}.
\end{equation}
Here $\langle\cdot, \cdot\rangle$ denotes the Euclidean inner product. Note that the evolution \eqref{eq:EFEvolution} is \emph{geometric}, i.e.\ if a smooth $f$ satisfies \eqref{eq:EFEvolution}, then for any smooth family of reparametrizations $\Phi\colon[0,T)\times I\to I$ so does $\hat{f}(t,x) \defeq (f\circ\Phi)(t,x)\defeq f(t, \Phi(t,x))$.
In addition to the evolution \eqref{eq:EFEvolution}, we prescribe \emph{clamped boundary conditions}, fixing position and the {unit tangent}
 of the curve at the endpoints of $I$. For an {immersed curve} 
 $f_0$ we hence study the following initial boundary value problem.
 \begin{align}
\label{eq:EF}
\left\{\begin{array}{rll}
\partial_tf &= - \nabla_s^2 \vKap -\frac{1}{2} \abs{\vKap} ^2 \vKap + \lambda \vKap + \theta\partial_s f & \text{ on } (0,T)\times  {I}  \\
f(0,x)&=f_0(x) & \text{ for }x\in {I} \\
f(t,y)&=p_y & \text{ for }0\leq t<T, y\in \partial I \\
\partial_s f(t,y) &= \tau_y &\text{ for } 0\leq t < T, y\in \partial I,\\
\end{array}
\right.
\end{align}
where the unknown $\theta \colon [0,T)\times I\to\R$, $\theta = \langle\partial_t f, \partial_s f\rangle$ is the \emph{tangential velocity}. By the integral representation of $\lambda$, \eqref{eq:EF} becomes a \emph{nonlocal quasilinear system} which is also \emph{degenerate} parabolic by its geometric nature. We assume that the boundary {data} $p_y\in \R^d, \tau_y\in \mathbb{S}^{d-1}\subset \R^d$ satisfy the \emph{compatibility conditions}
\begin{align}\label{eq:BCCompatibility}
	f_0(y) = p_y \text{ and } \partial_s f_0(y)=\tau_y \quad\text{ for }y\in \partial I.
\end{align}
Note that \eqref{eq:EF} is preserved under a smooth family of reparametrizations $\Phi$ which keeps the boundary  $\partial I$ fixed, {where the tangential velocity might change.} 

It is not difficult to see that $\lambda$ is chosen exactly in such a way that the length remains fixed {during the flow}, since along any sufficiently smooth solution of \eqref{eq:EF} we have
\begin{align}
	\label{eq:LengthConstantDuringFlow}
	\frac{\diff}{\diff t} \CalL(f) = -\int_I\langle\vKap, \partial_t f\rangle\diff s = \int_I \langle \nabla_s^2 \vKap + \frac{1}{2} \abs{\vKap} ^2 \vKap, \vKap\rangle \diff s - \lambda \int_I \abs{\vKap}^2\diff s =  0,
\end{align}
whereas the energy indeed decreases since by \eqref{eq:LengthConstantDuringFlow}
\begin{align}
	\label{eq:FlowMonotonical}
	\frac{\diff}{\diff t} \CalE(f)
	&= \int_I\langle \nabla \CalE(f), \partial_t f \rangle \diff s
	= \int_I\langle \nabla \CalE(f) - \lambda\vKap, \partial_t^{\perp} f \rangle \diff s = -\int_I \abs{\partial_t^{\perp} f}^2\diff s,
\end{align}
using {that the $L^2(\diff s_f)$-gradient of $\CalE$ is given by} $\nabla\CalE(f) = \nabla_s^2 \vKap + \frac{1}{2} \abs{\vKap} ^2 \vKap$. In the above calculations, we also used the fact that all boundary terms vanish due to the boundary conditions. In order for $\lambda$ to be well-defined, we need to ensure that $f(t)\defeq f(t, \cdot)$ is not a piece of a straight line. This can be guaranteed {with no restrictions on $\tau_0, \tau_1$} by requiring
\begin{align}\label{eq:f_0assumption}
\abs{p_0-p_1} < \ell \defeq \CalL(f_0),
\end{align}
{so $\CalE(f_0)>0$, see \Cref{subsec:lambda} for a more detailed analysis of $\lambda$.}

In \cite{DPL17}, long time existence for smooth solutions of \eqref{eq:EF} with tangential velocity $\theta\equiv 0$ under assumption \eqref{eq:f_0assumption} was shown with the help of interpolation inequalities. For the short time existence the authors of \cite{DPL17} refer to the beginning of Section 3 in \cite{DKS}, where the short time existence in the setting of H\"older spaces is only sketched for the case of closed curves.  Moreover, the uniform bounds in \cite[Theorem 1.1]{DPL17} imply subconvergence after reparametrization as $t\to\infty$. However, different sequences could still have different limits.

The contribution of this paper is twofold: First, we give a rigorous and fairly concise proof of short time existence and parabolic smoothing for the elastic flow \eqref{eq:EF}. 
Compared to the previous classical existence results for elastic flows, where the initial datum is assumed to be smooth \cite{PoldenPDE,DKS,DPL17} or at least with H\"older continuous second derivative \cite{KZE}, one major improvement is that we allow for rough initial values, lying merely in the natural \emph{energy space}, see \Cref{rem:p=2 justify} for a detailed discussion.  In contrast to existence theorems relying on the minimizing movement scheme (cf.\ \cite{Marius_Obstacle20,Okabe_Pozzi_Wheeler_p,Badal2022,Blatt_Hopper_Vorderobermeier22,Okabe_Wheeler_p}), our methods rely on \emph{maximal regularity}, \blue{yielding here smooth solutions}, cf.\ \Cref{thm:STE Main} below, while still allowing for rough initial data.
The price for this substantial improvement is that the necessary contraction estimates become quite technical and rely delicately on the precise structure of \eqref{eq:EF}. This is the first existence result for an elastic flow with general initial data of such weak regularity.

\begin{thm}\label{thm:STE Main}
	Let $f_0\in W^{2,2}(I;\R^d)$ be immersed, let ${p_0, p_1\in \R^d}$ and $\tau_0, \tau_1 \in \mathbb{S}^{d-1}$ satisfy \eqref{eq:BCCompatibility} and \eqref{eq:f_0assumption}. Then, there exists $T>0$ and a solution ${f\in W^{1,2}(0,T;L^2(I,\R^d))\cap L^2(0,T;W^{4,2}(I;\R^d))}$ of \eqref{eq:EF}.
\end{thm}
Moreover, we show that under the assumptions \eqref{eq:BCCompatibility} and \eqref{eq:f_0assumption}, the solution in \Cref{thm:STE Main} instantaneously becomes smooth, both in space and time, cf.\ \Cref{thm:Smoothing Main}.

Secondly we prove and apply a \emph{constrained {\L}ojasiewicz--{S}imon gradient inequality} (cf.\ \cite{Fabian}) to deduce convergence of the flow, where a new estimate (see \Cref{lem:BeschtesLemma}) 
substantially simplifies the argument for the convergence result compared to previous works, cf.\ \cite{CFS09, Loja}.

\begin{thm}\label{thm:Convergence Main}
	Let $f_0\in W^{2,2}(I;\R^d)$ be an immersed curve and suppose ${p_0, p_1\in \R^d}$ and $\tau_0, \tau_1 \in \mathbb{S}^{d-1}$ satisfy \eqref{eq:BCCompatibility} and \eqref{eq:f_0assumption}. Then, there exists a smooth family of curves $f\colon (0,\infty)\times I\to\R^d$ solving \eqref{eq:EF}, such that
		\begin{enumerate}[(i)]
		\item $f(t)\to f_0$ in $W^{2,2}(I;\R^d)$ as $t\to 0$;
		\item $f(t)\to f_{\infty}$ smoothly after reparametrization as $t\to\infty$, where $f_\infty$ is a  constrained clamped elastica, i.e.\ a solution of 
		\begin{equation}
		\left\{\begin{array}{rll}
		- \nabla_s^2 \vKap -\frac{1}{2} \abs{\vKap} ^2 \vKap + \lambda \vKap &= 0& \text{ on } I\\
		f(y)&=p_y &\text{ for }y\in \partial I \\
		\partial_s f(y) &= \tau_y &\text{ for } y\in \partial I
		\end{array}
		\right.
		\label{eq:ElasticaFixedLength}
		\end{equation}
		for some $\lambda \in \R$.
	\end{enumerate}
\end{thm}

{Together with the previously mentioned work \cite{DPL17} this paper completes the study of the existence and convergence of the elastic flow of clamped curves with fixed length. Unfortunately, due to the low regularity of the initial curves considered here, we are not able to show uniqueness for the solution of the geometric evolution equation \eqref{eq:EF}.}

In the smooth category, one can show uniqueness ``up to reparametrization'' by a PDE argument similar to \cite{Garcke2000}. {However, due to our low regularity we were not able to prove sufficient contraction estimates.} The reason for that is the rigid characterization of Lipschitz properties of Nemytskii operators, see for instance \cite[Theorem 3.10, Theorem 7.9]{Appell90}.

The elastic energy of curves has already been studied by Bernoulli. The analysis of the elastic flow, i.e.\ the one-dimensional analogue of the Willmore flow started with \cite{PoldenPDE} and \cite{DKS}. 
The boundary value problem for the elastic flow was considered in \cite{LIN20126414} for clamped curves and in \cite{dall2014evolution} for natural second-order boundary conditions, see also \cite{Wen,LinLueSchwetlick,Lin_Lue_2018} for related second order evolutions. For further related literature on elastic flows, we refer to \cite{NovagaOkabe,KZE,Okabe_Pozzi_Wheeler_p,Badal2022,Blatt_Hopper_Vorderobermeier22,Okabe_Wheeler_p}.
Recent research has also studied the geometric evolution
 of networks and previously achieved results were applied to the elastic flow of networks, see e.g.\ \cite{Elasticflowofnetworkslongtimeexistenceresult,GMP19,NovagaPozzi_p_network,GMP_network_LTE,DLP_network_STE}. Moreover, the elastic flow with different ambient geometries has been considered in \cite{Theelasticflowofcurvesonthesphere,MariusAdrian,Pozzetta}, especially, the case of hyperbolic space \cite{MariusAdrian} is of interest, cf.\  \cite{LangerSinger2,2005.13500}. Addtionally, we mention the elastic flow of closed curves under a length and area constraint \cite{Okabe_area}.

 The {\L}ojasiewicz--Simon gradient inequality is a remarkable result on \emph{(real) analytic} functions which was first proven in $\R^d$ \cite{Loja65}, and later generalized to infinite dimensions \cite{Simon83}, see also \cite{Chill03}. Nowadays, it is the fundamental tool for investigating the asymptotic properties of gradient flows with analytic energies, which has been used for many geometric evolution equations, {see for instance \cite{CFS09,Loja,Feehan,MantegazzaPozzetta,Rupp_Vol_pres,Rupp_Iso,Mantegazza_Pozzetta_Hypersurfaces,okabe2023convergence} and also \cite{NovagaOkabeCrelle} for a different approach.} The fixed-length constraint in \eqref{eq:EF} and \eqref{eq:LengthConstantDuringFlow} obstructs the use of \cite{Chill03} to deduce the gradient inequality, which is why we apply a recent extension to constrained energies \cite{Fabian}. We emphasize that this article is the first application of the constrained {\L}ojasiewicz--Simon gradient inequality for a constrained gradient flow.

This article is structured as follows. In \Cref{sec:STE}, we pick a specific tangential velocity such that \eqref{eq:EF} becomes a parabolic system, which we reduce to a fixed point equation. The existence of a fixed point is then established on a small time interval, using the concept of maximal $L^p$-regularity together with appropriate contraction estimates. \Cref{sec:Smoothing} is devoted to show instantaneous smoothing of our solution, both in space and time. After that, we prove long-time existence and a refined {\L}ojasiewicz--Simon gradient inequality to finally prove \Cref{thm:Convergence Main} in \Cref{sec:Convergence}.
For the sake of readability, some details on the contraction estimates and the parabolic smoothing have been moved to the appendix or can be found in the first author's dissertation \cite[Chapter 3]{Diss}.

\section{Short time existence}\label{sec:STE}
The goal of this section is to prove \Cref{thm:STE Main}. As in \cite{GMP19}, we prescribe an explicit tangential motion to transform \eqref{eq:EF} into a quasilinear parabolic {system}. We then perform a linearization and use the theory of maximal $L^p$-regularity and suitable contraction estimates to prove \Cref{thm:STE Main} using a fixed point argument. 
{We consider an initial datum merely lying in $W^{2,2}_{Imm}(I;\R^d)$, the space of $W^{2,2}$-immersions. This is a natural space for the elastic energy, since it is the roughest Sobolev space where $\CalE$ remains finite.}

\subsection{On the Lagrange multiplier}\label{subsec:lambda}
{To ensure that the Lagrange multiplier is well-defined, one needs to prevent the denominator from vanishing. \blue{Write} $\lambda(f) \eqdef \frac{N(f)}{2\CalE(f)}$, where $N(f)$ denotes the numerator in \eqref{eq:lambda} and observe that for a solution of \eqref{eq:EF} we have 
\begin{align*}
	\abs{f(t,0)-f(t,1)} = \abs{p_0-p_1} <\ell =\CalL(f(t)) \quad \text{ for all }t\in [0,T),
\end{align*}
using the boundary conditions, \eqref{eq:f_0assumption} and \eqref{eq:LengthConstantDuringFlow}. In particular, $f(t)$ cannot be part of a straight line, so $\CalE(f(t))>0$ for all $t\in [0,T)$. Moreover, we observe that after integration by parts we have
\begin{align}\label{eq:lambdaIBP}
N(f) =\int_I\langle \nabla \mathcal{E}(f), \vKap\rangle\diff s 
=
\langle \nabla_s \vKap, \vKap \rangle |_{\partial I} - \int_I
|\nabla_s \vKap|^2 \diff s + \frac{1}{2} \int_I |\vKap|^4 \diff s.
\end{align}
Note that in \eqref{eq:lambdaIBP}, no derivatives of second order of the curvature appear, which means that the Lagrange multiplier is formally of lower order compared to $\nabla\CalE(f)$. This \blue{is} extremely useful later on, since we can rely on the well-studied property of maximal $L^p$-regularity for a local operator in the linearization and treat the Lagrange multiplier as a nonlinearity in the fixed point argument.
}
\subsection{From the geometric problem to a quasilinear PDE}\label{subsec:quasilin}
 As a next step, we explicitly compute the right hand side of \eqref{eq:EFEvolution}. By \Cref{prop:ZeugsInKoordinaten}
\begin{align*}
\nabla\CalE(f) &=  \nabla_s^2 \vKap +\frac{1}{2} \abs{\vKap} ^2 \vKap =\CalA(f)^{\perp},
\end{align*}
where
\begin{align}\label{eq:defA}
\nonumber
\mathcal{A}(f) &\defeq  \frac{\partial_x^{4}f}{\gamma^{4}} -6 \frac{\langle\partial_{x}^{2} f, \partial_x f\rangle}{\gamma^{6}}\partial_x^{3}f - 4 \frac{\langle\partial_x^{3}f, \partial_x f\rangle}{\gamma^{6}}\partial_x^{2}f - \frac{5}{2} \frac{\abs{\partial_x^{2} f}^{2}}{\gamma^{6}} \partial_x^{2} f \\
&\quad + \frac{35}{2} \frac{\langle\partial_x^{2}f, \partial_x f\rangle^{2}}{\gamma^{8}} \partial_x^{2}f\nonumber\\
& \eqdef  \frac{\partial_x^{4}f}{\gamma^{4}} +\tilde{F}(\gamma^{-1}, \partial_x f, \partial_x^2 f, \partial_x^3 f).
\end{align}
In order to solve \eqref{eq:EF}, we study the following evolution problem, prescribing an explicit tangential motion $\theta =\mu$ to make the problem parabolic. {We want to find a family of immersions $f\colon [0,T)\times I\to\R^d$ satisfying}
\begin{align}\label{eq:PludaDeTurck}
\left\{\begin{array}{rll}
\partial_t f &= - \nabla_s^2 \vKap -\frac{1}{2} \abs{\vKap} ^2 \vKap + \mu \partial_s f + \lambda \vKap & \text{ on } (0,T)\times {I}  \\
f(0,x)&=f_0(x) & \text{ for }x\in {I} \\
f(t,y)&=p_y &  \text{ for }0\leq t<T,y\in \partial I \\
\partial_{x} f(t,y)&= \tau_y{\gamma_0(y)} &  \text{ for }0\leq t < T,y\in \partial I.\\
\end{array}\right.
\end{align} 
with $\lambda$ as in \eqref{eq:lambda} and $\mu=\mu(f)\colon [0,T)\times I\to \R$ given by $\mu\defeq - \langle\CalA(f), \partial_s f\rangle$. Note that the first order boundary conditions are a linear version of the general boundary conditions in \eqref{eq:EF}, and thus easier to handle. The system \eqref{eq:PludaDeTurck} is often referred to as the \emph{analytic problem}.

For $1<p<\infty$ and $T>0$, we consider the space of solutions
\begin{align*}
	\X_{T,p} \defeq W^{1,p}\left(0,T; L^p(I;\R^d)\right)\cap L^p\left(0,T; W^{4,p}(I;\R^d)\right)
\end{align*}
and the space of data
\begin{align*}
	\Y_{T,p}^1 \defeq L^p\left(0,T;L^p(I;\R^d)\right).
\end{align*}
The space of initial data is given by the \emph{Besov space}
\begin{align*}
	\Y_p^{2}\defeq \{ f(0) \mid f\in \X_{T,p}\} = B^{4(1-\frac{1}{p})}_{p,p}(I;\R^d),
\end{align*}
see for instance {\cite[Section 2]{DHP07}}. We also consider the solution space with vanishing trace at time $t=0$ given by
\begin{align*}
	\prescript{}{0}{\X}_{T,p} \defeq \{ f\in \X_{T,p}\mid f(0)=0\}.
\end{align*}
For convenience, we also set $\Y_{T,p} \defeq \Y_{T,p}^1 \times \Y_p^2$.

\subsection{Linearization of the analytic problem}
If we linearize \eqref{eq:PludaDeTurck} for $\lambda\equiv 0$, we obtain a linear parabolic system. This system is a local PDE which we can apply maximal regularity theory to. First, assuming $\lambda\equiv 0$ and using \eqref{eq:defA}, the evolution in \eqref{eq:PludaDeTurck} has the form
\begin{align*}
	\partial_t f  = -\CalA(f) \eqdef -\frac{\partial_x^{4}f}{\gamma^{4}} - \tilde{F}(\gamma^{-1},\partial_xf, \partial_{x}^{2}f, \partial_x^{3}f)
\end{align*}
with $\mathcal{A}$ as in \eqref{eq:defA}. If we freeze coefficients for the highest order term at the initial datum $f_0$ we get
 \begin{align}\label{eq:linearization}
 	\partial_t f + \frac{\partial_x^{4}f}{{\gamma_0}^{4}} &= \left(\frac{1}{\gamma_0^{4}}-\frac{1}{\gamma^{4}}\right)\partial_x^{4} f- \tilde{F}(\gamma^{-1},\partial_xf, \partial_{x}^{2}f, \partial_x^{3}f)\nonumber \\
 	&\eqdef F(\gamma^{-1},\partial_x f, \partial_x^{2}f, \partial_x^{3}f, \partial_x^{4}f),
 \end{align}
 where $\gamma_0 \defeq \gamma(0, \cdot) = \abs{\partial_x f_0}$ and $\tilde{F}$ \blue{is} as in \eqref{eq:defA}. The linearized system we associate to \eqref{eq:PludaDeTurck} with $\lambda\equiv 0$ is 
 \begin{align}\label{eq:linearSystem}
 	\left\{\begin{array}{rll}
 	\partial_t f + \frac{1}{\gamma_0^{4}}\partial_x^{4}f &= F& \text{ on } (0,T)\times {I}  \\
 	f(0,x)&=f_0(x) & \text{ for }x\in {I} \\
 	f(t,y)&=p_y &  \text{ for }0\leq t<T,y\in \partial I \\
	{\partial_x f(t,y)}&= \tau_y{\gamma_0(y)} &  \text{ for }0\leq t < T,y\in \partial I.\\
 	\end{array}\right.
 \end{align}

 We can now apply the general $L^p$-theory for parabolic systems to obtain the following classical maximal regularity result, whose proof can be found in \cite[Chapter 3, Section 2.3]{Diss}. For the definition of the {spaces for the boundary data}  $\CalD^{i}_{T,p}$ with $i=0,1,$ see \eqref{eq:BoundarySpaces}. 
 
 \begin{thm}\label{thm:maxReg} 
 	Let $p\in (\frac{5}{3},\infty)$, $0<T\leq T_0$. Suppose $a\in \CalC([0,T_0]\times I;\R)$ such that $a(t,x)\geq \alpha$ for some $\alpha>0$ and all $t\in [0,T_0], x\in I$.
 	Let $(\psi,f_0) \in \Y_{T,p}$, $b^0 \in \CalD^0_{T,p}$ and $b^1\in \CalD^1_{T,p}$ such that the following compatibility conditions are satisfied:
 	\begin{align}\label{eq:Compatibility}
 	b^0(0,y) &= f_0(y) &\text{for }y\in \partial I,\nonumber \\
 	b^1(0,y) &= \partial_x f_0(y) &\text{for }y\in \partial I.
 	\end{align}
 	Then, there exists a unique $f\in \X_{T,p}$ such that
 	\begin{align}\label{eq:LinPDESys}
 	\left\{\begin{array}{rll}
 	\partial_t f + a\partial_x^{4}f &= \psi& \text{ on } (0,T)\times {I}  \\
 	f(0,x)&=f_0(x) &  \text{ for }x\in {I} \\
 	f(t,y)&=b^0(t,y) &  \text{ for }0\leq t<T,y\in \partial I \\
 	\partial_x f(t,y)&= b^1(t,y) &  \text{ for }0\leq t < T,y\in \partial I,\\
 	\end{array}\right.
 	\end{align}
	 and there exists $C=C(p, T, a)>0$ such that
 	\begin{align}\label{eq:LinearAPrioriEstim}
 	\norm{f}{\X_{T,p}} \leq C\left( \norm{\psi}{\Y_{T,p}^1} + \norm{f_0}{\Y_{p}^2} + \norm{b^0}{\CalD^0_{T,p}} + \norm{b^1}{\CalD^1_{T,p}}\right).
 	\end{align}
 	{Moreover, if $b^0=0$ and $b^1=0$, then we may choose $C=C(p,T_0,a)$ independent of $T\leq T_0$.}
 \end{thm}
 
Now, we want to solve \eqref{eq:PludaDeTurck} for initial data $f_0\in W^{2,2}(I;\R^d)$ using a fixed point argument. Note that $B^{2}_{2,2}(I;\R^d)=W^{2,2}(I;\R^d)$ by \eqref{eq:Besov=Sobolev}, so $p=2$ is a fine setup to deal with the desired initial data, see \Cref{rem:p=2 justify} for a more detailed discussion.
We observe that the linearized system \eqref{eq:linearSystem} can be viewed as a special case of \Cref{thm:maxReg} with $a= \frac{1}{\gamma_0^4}$, $b^0 = (p_0, p_1)$, $b^1 = (\tau_0, \tau_1)$ and $\psi=F$. 

 Throughout the rest of this section, we exclusively work with $p= 2$. {To simplify notation} the spaces $\X_{T}, \Y_T, \CalD^0, \CalD^1$ denote the respective spaces with $p=2$. 
 
\subsection{Contraction estimates}\label{subsec:Contraction} 
The key ingredient in the proof of the short time existence is a contraction estimate for the nonlinearity in \eqref{eq:PludaDeTurck}. We fix an initial datum $f_0\in W^{2,2}_{Imm}(I;\R^d)$ and boundary conditions $p_0, p_1\in \R^d$ and $\tau_0, \tau_1\in \mathbb{S}^{d-1}$ satisfying  \eqref{eq:BCCompatibility} and \eqref{eq:f_0assumption}. 
\blue{For a \emph{reference flow} $\bar{f}\in \X_{T=1}$ with $\bar{f}(0)=f_0$, and some $M$ and $T\in (0,1]$ we define}
\begin{align}\label{eq:def B_M}
\blue{\bar{B}_{T,M}\defeq \left\{f \in \X_T \mid f(0)=f_0 \text{ and } \norm{f-\bar{f}}{\X_T}\leq M\right\}.}
\end{align}
We denote by $T$ the \emph{existence time} and by $M$ the \emph{contraction radius.}  \blue{Since we take $T,M>0$ small later on, it is no restriction to only consider $T,M\leq 1$.} 
Later, we choose a specific reference flow $\bar{f}$, see \Cref{def:referenceSol}.

First, the following lemma yields uniform bounds from below on the arc-length element and the elastic energy for small times, \blue{ensuring that the system \eqref{eq:PludaDeTurck} does not immediately become singular.} A detailed proof can be found in \cite[Chapter 3, Section 2.4.1]{Diss}.


\begin{lem}\label{lem:boundsgamma}
	For \blue{$T=T(\bar{f})\in (0,1]$} small enough and \blue{$M\in (0, 1]$,} any $f\in \blue{\bar{B}_{T,M}}$ satisfies	 $\gamma(t,x)\geq \inf_{I}\frac{\gamma_0}{2}$ for all $(t,x)\in [0,T)\times I$.	In particular, all curves $f(t,\cdot)$ are immersed.
\end{lem}
\begin{lem}\label{lem:lowerBoundsEnergy}
	\blue{For $T=T(\bar{f})\in (0,1]$ small enough and $M\in (0,1]$, any $f\in \blue{\bar{B}_{T,M}}$ satisfies} $\CalE(f(t)) \geq \frac{\CalE(f_0)}{3}>0$ (cf.\ \eqref{eq:f_0assumption}) for all $t\in [0,T)$, .
\end{lem}


We now state the crucial contraction property of the nonlinearities. Since the space of initial data is the energy space, cf.\ \Cref{rem:p=2 justify}, the necessary estimates are quite involved and rely on the special structure of \eqref{eq:EFEvolution}. For the sake of readability, some of the details and the proof of the following lemma are moved to \Cref{app:contr}.

\begin{lem}\label{lem:NContraction} Let $q\in (0,1)$. Then the following maps 
	\begin{align*}
	\begin{split}
	\CalF \colon \blue{\bar{B}_{T,M}} \to \Y_T^1,\quad \CalF(f) &\defeq F(\gamma^{-1}, \partial_x f, \partial_x^2 f, \partial_x^3 f, \partial_x^4f) \\
	\Lambda \colon \blue{\bar{B}_{T,M}} \to \Y_T^1,\quad \Lambda(f) &\defeq \lambda(f)\vKap_f\\
	\mathcal{N}\colon \blue{\bar{B}_{T,M}} \to \Y_T^1, \quad \mathcal{N}(f)  &\defeq \CalF(f)+\Lambda(f),
	\end{split}
	\end{align*}
	are well-defined $q$-contractions (i.e.\ Lipschitz \blue{continuous} with Lipschitz constant $q$) for \blue{$T=T(q,\bar{f})$,} \blue{$M=M(q, \bar{f})\in (0,1]$} small enough,
	with $F$ as in \eqref{eq:linearization} and $\lambda$ as in \eqref{eq:lambda}.
\end{lem}
\subsection{The fixed point argument}

We now reduce the analytic problem \eqref{eq:PludaDeTurck} to a fixed point equation and show local existence and uniqueness via the contraction principle. To that end, we first choose a specific {reference solution} $\bar{f}$ in \eqref{eq:def B_M} on  \blue{the time interval $[0,1]\supset [0,T]$ for $0<T\leq 1$.}
\begin{defi}\label{def:referenceSol}
	 We define the \emph{reference solution} $\bar{f}$ to be the unique solution of the following initial boundary value problem.
	\begin{align*}
		\left\{\begin{array}{rll}
		\partial_t \bar{f} + \frac{\partial_x^{4}\bar{f}}{\gamma_0^4} &= 0 & \text{ on } \blue{[0,1)}\times {I}  \\
		\bar{f}(0,x)&=f_0(x) & \text{ for }x\in {I} \\
		\bar{f}(t,y)&=p_y & \text{ for }0\leq t<\blue{1}, y\in \partial I \\
		\partial_x \bar{f}(t,y)&= \tau_y {\gamma_0(y)} &\text{ for } 0\leq t < \blue{1},y\in \partial I.\\
		\end{array}\right.
	\end{align*}
\end{defi}

	\blue{Existence and uniqueness in the class $$W^{1,2}\left(0,\blue{1}; L^2(I;\R^d)\right)\cap L^2\left(0,\blue{1}; W^{4,2}(I;\R^d)\right)$$ follows from \Cref{thm:maxReg}. {Note that the restriction of the solution to any time interval $[0,T]$ is the unique solution in the class $\X_{T}$ for all $0<T\leq \blue{1}$}.
	}

	 \blue{Fix $q\in (0,1)$ and take ${T}={T}(q,\bar{f})\in (0,1],M=M(q,\bar{f})\in (0,1]$} small enough such that \Cref{lem:boundsgamma,lem:lowerBoundsEnergy,lem:NContraction} hold. Let $f\in \blue{\bar{B}_{T,M}}$. Then, we have $\mathcal{N}(f) \in \Y_{{T}}^{1}$, cf.\ \Cref{lem:NContraction}. For $\psi\defeq \mathcal{N}(f), b^0 \defeq (p_0,p_1)$, $b^1 \defeq (\tau_0, \tau_1)$, $a\defeq \gamma_0^{-4}\in \CalC([0,\blue{1}]\times I)$ 
	the compatibility conditions \eqref{eq:Compatibility} are satisfied, since by \eqref{eq:BCCompatibility} we have
	\begin{align*}
	b^0(0,y) &= f_0(y)&\text{for }y\in \partial I, \\
	b^1(0,y) &= \tau_y {\gamma_0(y)} = \partial_x f_0(y)&\text{for }y\in \partial I.
	\end{align*}
	Hence, by \Cref{thm:maxReg}, there exists a unique solution $g\in \X_T$ of the linear initial boundary value problem
	\begin{align}\label{eq:FixedPointDefi}
	\left\{\begin{array}{rll}
	\partial_t g + \frac{\partial_x^{4}g}{\gamma_0^4} &= \mathcal{N}(f) & \text{ on } (0,T)\times {I}  \\
	g(0,x)&=f_0(x) & \text{ for }x\in {I} \\
	g(t,y)&=p_y & \text{ for } 0\leq t<T,y\in \partial I \\
	\partial_x g(t,y) &= \tau_y {\gamma_0(y)} & \text{ for }0\leq t < T,y\in \partial I.\\
	\end{array}\right.
	\end{align}
\begin{defi}\label{def:FixedPointMap}
	We define the map $\Phi\colon \blue{\bar{B}_{T,M}} \to \X_T, \Phi(f)\defeq g$, where $g\in \X_T$ is the unique solution to \eqref{eq:FixedPointDefi}.
\end{defi}

\begin{remark}
	Finding a solution of \eqref{eq:PludaDeTurck} in the ball $\blue{\bar{B}_{T,M}} \subset \X_T$ is equivalent to finding a fixed point of the map $\Phi$ in \Cref{def:FixedPointMap}.
\end{remark}

We now show that $\Phi$ is a contraction on $\blue{\bar{B}_{T,M}}$ for $T,M>0$ small enough.

\begin{prop}\label{prop:PhiContraction}
	Let $q\in (0,1)$. Then there exist \blue{$M=M(q,\bar{f})\in (0,1]$, $T=T(q,M,\bar{f})\in (0,1]$} such that $\Phi \colon\blue{\bar{B}_{T,M}} \to \blue{\bar{B}_{T,M}}$ is well-defined and a $q$-contraction, i.e.
	\begin{align}\label{eq:PhiContraction}
		\norm{\Phi(f)-\Phi(\tilde{f})}{\X_T}\leq q \norm{f-\tilde{f}}{\X_T}
	\end{align}
	for all $f,\tilde{f} \in \blue{\bar{B}_{T,M}}$.
\end{prop}
\begin{proof}
	\blue{
	\underline{The contraction property:} Let $q\in (0,1)$ and $f, \tilde{f}\in \blue{\bar{B}_{T,M}}$ and let $g=\Phi(f)$, $\tilde{g}=\Phi(\tilde{f})$. We observe that $g-\tilde{g}$ vanishes at \blue{$t=0$ and at} the boundary $\partial I$ up to first order. Hence, by \Cref{def:FixedPointMap} and \eqref{eq:LinearAPrioriEstim}, \blue{for some $C=C(f_0)=C(\bar{f})>0$, independent of $T\in (0,1]$, we have}
	\begin{align}\label{eq:PhiContraction1}
		\norm{g-\tilde{g}}{\X_T}&\leq C \norm{\mathcal{N}(f)-\mathcal{N}(\tilde{f})}{L^2(0,T;L^{2})}.
	\end{align}
	Taking $T=T(q,\bar{f}),M=M(q,\bar{f})\in (0,1]$ small enough so that \Cref{lem:NContraction} can be applied with $q$ replaced by $q/(2C)$, we have
	\begin{align}\label{eq:PhiContraction2}
		\norm{\mathcal{N}(f)-\mathcal{N}(\tilde{f})}{L^2(0,T;L^2)}\leq \frac{q}{2C} \norm{f- \tilde{f}}{\X_T}.
	\end{align}
	 Equations \eqref{eq:PhiContraction1} and \eqref{eq:PhiContraction2} imply \eqref{eq:PhiContraction}.}
	
	\underline{Well-definedness:} 
	Let $f\in \blue{\bar{B}_{T,M}, g=\Phi(f)}$. \blue{Again by  \eqref{eq:LinearAPrioriEstim} we find}
		\begin{align}\label{eq:PhiWellDef}
		\norm{g-\bar{f}}{\X_T} &\leq C \norm{\mathcal{N}(f)-0}{L^2(0,T;L^{2})}\nonumber \\
		&\leq C \left(\norm{\mathcal{N}(f)-\mathcal{N}(\bar{f})}{L^2(0,T;L^2)} + \norm{\mathcal{N}(\bar{f})}{L^2(0,T;L^{2})} \right) \nonumber\\
		&\leq \frac{q}{2}\norm{f-\bar{f}}{\X_T} + C \norm{\mathcal{N}(\bar{f})}{L^2(0,T;L^{2})}\nonumber\\
		&\leq \frac{q}{2}M + C\norm{\mathcal{N}(\bar{f})}{L^2(0,T;L^{2})},
	\end{align}
	\blue{where we applied \eqref{eq:PhiContraction2} with $\tilde{f}=\bar{f}$ in the third step.} Now, by dominated convergence we have $\norm{\mathcal{N}(\bar{f})}{L^2(0,T;L^2)} \leq \frac{M}{2C}$ reducing \blue{$T=T(q,M, \bar{f})\in (0,1]$} if necessary.
%
Then, from \eqref{eq:PhiWellDef} we conclude $\norm{\Phi(f)-\bar{f}}{\X_T}\leq M$.
\end{proof}

\begin{thm}\label{thm:STE}
	{Let $f_0\in W^{2,2}_{Imm}(I;\R^{d})$, $p_0, p_1\in \R^d$, $\tau_0, \tau_1\in \mathbb{S}^{d-1}$ satisfying \eqref{eq:BCCompatibility} and \eqref{eq:f_0assumption}.} Then there exists $M>0$ and $T>0$ such that the system \eqref{eq:PludaDeTurck} has a unique solution $f\in \blue{\bar{B}_{T,M}}\subset  W^{1,2}\left(0,T; L^2(I;\R^d)\right)\cap L^2\left(0,T; W^{4,2}(I;\R^d)\right)$.
\end{thm}
\begin{proof}
For $T,M>0$ as in \Cref{prop:PhiContraction} with $q=\frac{1}{2}$, the map $\Phi\colon \blue{\bar{B}_{T,M}}\to\blue{\bar{B}_{T,M}}$ is a contraction in the complete metric space $\blue{\bar{B}_{T,M}}$ and hence has a unique fixed point $f\in \blue{\bar{B}_{T,M}}$ by the contraction principle. Since any fixed point of $\Phi$ is a solution of \eqref{eq:PludaDeTurck} in $\blue{\bar{B}_{T,M}}$  and vice versa, the claim follows.
\end{proof}

\begin{remark}\label{rem:STESolBounds}
	By the construction of our solution and \Cref{lem:boundsgamma} and \Cref{lem:lowerBoundsEnergy} the arc-length element $\abs{\partial_x f}$ and the elastic energy of the solution $f$ in \Cref{thm:STE} are bounded from below and above, uniformly in $t\in [0,T)$.
\end{remark}
This immediately implies \Cref{thm:STE Main}.
\begin{proof}[{Proof of \Cref{thm:STE Main}}]
	By \Cref{thm:STE}, there exists $T>0$ and  a solution $f$ of \eqref{eq:PludaDeTurck} such that $f\in W^{1,2}(0,T;L^2(I;\R^d))\cap L^2(0,T;W^{4,2}(I;\R^d))$. Consequently, $f$ solves \eqref{eq:EF}, since at the boundary we have
	\begin{align*}
		\partial_{s_f}f(t,y) = \frac{\partial_x f(t,y)}{\abs{\partial_x f(t,y)}} = \frac{\gamma_0(y)\tau_y}{\abs{\gamma_0(y)\tau_y}} = \tau_y \quad \text{ for }t\in[0,T), y\in \partial I.&\qedhere
	\end{align*}
\end{proof}
\begin{remark}\label{rem:p=2 justify}
	Our assumption on the regularity of the initial datum is very natural. On the one hand, the space $W^{2,2}_{Imm}(I;\R^d)$ is the correct energy space associated to the elastic energy, so we would like to obtain short time existence for an initial datum in $W^{2,2}_{Imm}(I;\R^d)$. In view of the linear problem in \Cref{thm:maxReg}, working in the Sobolev scale one would hence  {need} to pick $p\in (1,\infty)$ such that $W^{2,2}(I;\R^d)\hookrightarrow B^{4(1-\frac{1}{p})_{p,p}}(I;\R^d) =\Y_p^2$.
	
	However, in order to estimate the denominator of the Lagrange multiplier $\lambda$, we want continuity of our solution with values in $W^{2,2}(I;\R^d)$. Using \Cref{prop:SobolevEmbeddings} (i), this can {be achieved} if $B^{4(1-\frac{1}{p})_{p,p}}(I;\R^d) \hookrightarrow W^{2,2}(I;\R^d)$.
	
	Clearly, this can only work for $p=2$.  Moreover, for the same reason as above, {even the introduction of} time-weighted Sobolev spaces {would} not provide {solutions with lower initial regularity.}
\end{remark}

\begin{thm}
\label{thm:globalUnique}
The solution $f \in \blue{\bar{B}_{T,M}}
$ in \Cref{thm:STE} is the unique solution of \eqref{eq:PludaDeTurck} in the whole space $W^{1,2}\left(0,T; L^2(I;\R^d)\right)\cap L^2\left(0,T; W^{4,2}(I;\R^d)\right)$. 
\end{thm}

\begin{proof}
First we note that any restriction of the solution $f \in \bar B_M $ to a smaller time interval $[0,\tilde T]$ is again the unique solution of \eqref{eq:PludaDeTurck} in $\bar B_M$ on $[0,\tilde T]$ by \autoref{thm:STE}. Now, we let $T_1, T_2 > 0$ and assume that $f_i \in W^{1,2}\left(0,T_i; L^2(I;\R^d)\right)\cap L^2\left(0,T_i; W^{4,2}(I;\R^d)\right)$, $i = 1,2$ are two {families of immersions satisfying} \eqref{eq:PludaDeTurck} with $f_0 \in W^{2,2}_{{Imm}}(I)$. Without loss of generality we may assume that $T_1 \leq T_2$. We claim that $f_2 |_{[0,T_1]} = f_1$.\\
To show the claim we define $\bar{t} = \sup\{t \in [0,T_1): f_1(s) = f_2(s) \, \forall\,  0 \leq s \leq t\}$. Note that $\bar{t}$ is well-defined by \autoref{prop:SobolevEmbeddings} (i). We need to show that $\bar{t} = T_1$. To do so we first prove that $\bar{t} > 0$. Indeed, for $T \searrow 0$, we have $\norm{f_i|_{[0,T]}}{\X_T} \to 0$ by the dominated convergence theorem, and the same holds for the reference flow $\overline f$ from \Cref{def:referenceSol}. Thus, for $T  > 0$ small enough, $f_i|_{[0,T]} \in \bar B_M$ for $i = 1,2$. Further decreasing $T>0$ if necessary we obtain from \autoref{thm:STE} that $f_1|_{[0,T]} = f_2|_{[0,T]}$ is the unique solution $f \in \bar B_M$. Thus, $f_1(s) = f(s) = f_2(s)$ for all $0 \leq s \leq T$, showing that $\bar{t} \geq T > 0$.\\
We now assume that $\bar{t} < T_1$. Since $f_i \in \X_{T_1} \hookrightarrow BUC([0,T_1],W^{2,2}(I;\R^d))$ and both solutions are immersed for all times, we find that $f_0 \defeq f_1(\bar{t}) \in W^{2,2}_{Imm}(I;\R^d)$. Whence, by \autoref{thm:STE}, there exist $M >0$, $T > 0$ such that \eqref{eq:PludaDeTurck} has a unique solution $f \in \bar B_M$. Observing that $f_i(\bar{t} + t, \cdot)|_{0 \leq t \leq T_1-\bar{t}}$, $i = 1,2$ are both solutions to \eqref{eq:PludaDeTurck} with the same initial value $f_0$, we find by similar arguments as above that $f_1(\bar{t} + \cdot) = f = f_2(\bar{t} + \cdot) $ on $[0, T)$, contradicting the definition of $\bar{t}$.
\end{proof}

\section{Parabolic smoothing}
\label{sec:Smoothing}
The goal of this section is to show that our solution $f$ from \Cref{thm:STE} instantaneously becomes smooth. 
\begin{thm}\label{thm:Smoothing Main}
		 Let {$f_0\in W^{2,2}_{Imm}(I;\R^d)$ such that  \eqref{eq:BCCompatibility} and \eqref{eq:f_0assumption}} are satisfied. Then, there exists $0<T_1\leq T$ such that the solution $f$ in \Cref{thm:STE Main} is smooth on $(0,T_1)$, i.e.\ $f\in \CalC^{\infty}((0,T_1)\times I;\R^d)$.
\end{thm}
A close examination of the contraction estimates in \Cref{app:contr} reveals that the critical embeddings are used, for instance in \eqref{eq:dx^3L4L2}, \eqref{eq:trdx^2L8} and \eqref{eq:trdx^3L8/3}. Thus, higher integrability of the nonlinearity cannot be obtained by standard estimates relying on Hölder's inequality. Therefore, we cannot directly start the usual bootstrap argument to show smoothness. Instead, we use an instantaneous gain of regularity in the time variable, relying on Angenent's parameter trick \cite{Angenent1,Angenent2}, \blue{see also \cite{EPS} and \cite[Chapter 9]{MR3524106}. Since we want to conclude smoothness in space \emph{up to the boundary}, we first show increased regularity in time before deducing global smoothness up to the boundary by using parabolic Schauder theory.} For the sake of readability, we omit the proof of the following proposition and refer to \cite[Chapter 3, Section 3.1]{Diss}.

\begin{prop}\label{thm: f Cw in Time}
	Let $f\in W^{1,2}(0,T;L^2(I;\R^d))\cap L^2(0,T;W^{4,2}(I;\R^d))$ be the unique solution of \eqref{eq:PludaDeTurck}, given by \Cref{thm:STE}. Then there exists $0<T_1<T$ such that $f\in \CalC^{\omega}((0,T_1); W^{2,2}(I,\R^d))$.
\end{prop}

Next, we use
the higher time regularity to improve the integrability of $\lambda$, which \blue{ then allows} us to start a bootstrap argument.
First, we recall the following modification of  \cite[Lemma 4.3]{DPL17}.

\begin{lem}\label{lem:Lambda Trick Anna}
	Let $f\in \X_{T,2}$ be a solution of \eqref{eq:EF}. Then, we have
	\begin{align*}
		&\abs{\lambda}\big(\ell- \abs{p_1-p_0}\big)\leq 2\ell\norm{\partial_t^{\perp}f}{L^1(\diff s)}+\int_I\abs{\vKap}^2\diff s+\int_I\abs{\nabla_s\vKap}\diff s.
	\end{align*}
\end{lem}
\begin{proof}
	We proceed as in \cite[Lemma 4.3]{DPL17}.
	Let $l\colon [0,T)\times I\to \R^d$ be the parametrization of the line segment from $p_0$ to $p_1$ given by
	\begin{align*}
	l(t,x)\defeq p_0 + \frac{\varphi(t,x)}{\ell}(p_1-p_0),	
	\end{align*}
	with $\varphi(t, \xi)\defeq \int_{0}^\xi \abs{\partial_x f}\diff x$ for $(t,\xi)\in [0,T)\times I$. Then $l(t,0)=p_0$, $l(t,1)=p_1$ and $\partial_{s } l(t, \cdot) = \frac{1}{\ell}(p_1-p_0)$. Therefore, using $\nabla_s^2\vKap + \frac{1}{2}\abs{\vKap}^2\vKap = \nabla_s\left(\nabla_s\vKap + \frac{1}{2}\abs{\vKap}^2\partial_s f\right)$ (cf.\ \cite[p. 1048]{DPL17}), we find after integrating by parts
	\begin{align*}
		\int_{I} \langle\partial_t^{\perp} f, f-l\rangle \diff s  &= \left.\langle(\lambda-\frac{1}{2}\abs{\vKap}^2)\partial_s f - \nabla_s\vKap, f-l\rangle\right\vert_{\partial I} + \frac{1}{2} \int_I\abs{\vKap}^2\diff s  - \lambda \int_I\diff s  \\
		&\quad - \frac{1}{\ell} \int_{I}  \langle \nabla_s\vKap+\frac{1}{2}\abs{\vKap}^2\partial_s f - \lambda \partial_s f, p_1-p_0\rangle\diff s .
	\end{align*}
	Consequently, since $f=l$ on the boundary, we have
	\begin{align*}
		&\abs{\lambda}(\ell-\abs{p_1-p_0}) =  \left(1-\frac{\abs{p_1-p_0}}{\ell}\right)\abs{\lambda}\int_{I}\diff s
		\\
		&\leq \int_{I} \abs{\partial_t^{\perp}f}\diff s  \norm{f-l}{\infty} + \frac{1}{2}\int_{I}\abs{\vKap }^2\diff s  + \frac{\abs{p_1-p_0}}{2\ell}\int_I\abs{\vKap }^2\diff s   + \frac{\abs{p_1-p_0}}{\ell}\int_I\abs{\nabla_s\vKap }\diff s.
	\end{align*}
	Using \eqref{eq:f_0assumption} and the simple estimate $\norm{f-l}{\infty} \leq 2\ell$  yields the claim.
\end{proof}
Note that a priori, the Lagrange multiplier $\lambda$ \blue{is} only $L^2(0,T)$ for $f\in \X_{T,2}$. The next proposition improves this integrability, at least on a small timescale bounded away from zero.
\begin{lem}\label{prop:lambdaL^4}
	Let $f$ be the solution of \eqref{eq:PludaDeTurck} from \Cref{thm:STE} and let $T_1>0$ as in \Cref{thm: f Cw in Time}. Then, for any $0<\varepsilon<T_1$ we have $\lambda(f)\in L^4(\varepsilon,T_1)$.
\end{lem}
\begin{proof}
	As a consequence of \Cref{thm: f Cw in Time}, we have $\partial_t f\in \CalC^{\omega}((0,T_1);W^{2,2}(I;\R^d))$ and thus we get $\partial_t f\in \CalC^{0}([\varepsilon,T_1]\times I;\R^d)$. 
	Hence, \Cref{lem:Lambda Trick Anna} and \eqref{eq:f_0assumption} yield that $\lambda$ has the same integrability on $(\varepsilon, T_1)$ as $\int_I \abs{\nabla_s\vKap}\diff s$. By \Cref{prop:ZeugsInKoordinaten} (ii) and the uniform bounds on the arc-length element, cf.\ \Cref{rem:STESolBounds}, it suffices to show  $\partial_x^3 f \in L^4(\varepsilon,T_1;L^1(I;\R^d))$, since $\partial_x^2f\in \CalC^{0}([\varepsilon,T_1];L^2(I;\R^d))$. {In fact using \Cref{prop:SobolevEmbeddings2} (i) as in \eqref{eq:dx^3L4L2}, we even get $\partial_x^3 f\in L^4(\varepsilon, T_1; L^2(I;\R^d))$.}
	\end{proof}

	The improved integrability of $\lambda$ in \Cref{prop:lambdaL^4} enables us to start a bootstrap argument to increase the Sobolev regularity of our solution in \Cref{thm:STE}. Note that by Sobolev embeddings, in order to prove smoothness of our solution it suffices to reach $\X_{T,p}$ with $p>5$, see \Cref{lem:p5smooth}.
	
	\begin{lem}\label{prop:finMR20}
		Let $f$ be as in \Cref{thm:STE}, let $T_1>0$ be as in \Cref{thm: f Cw in Time} and let $0<\varepsilon<T_1$. Then $f\in W^{1,20}(\varepsilon, T_1; L^{20}(I;\R^d))\cap L^{20}(\varepsilon
		,T_1;W^{4,20}(I;\R^d))$ .
	\end{lem}
	\begin{proof}
		See \cite[Chapter 3, Section 3.3]{Diss}.
	\end{proof}

Finally, \Cref{thm:Smoothing Main} follows from parabolic \blue{Schauder} theory and the following
\begin{lem}
	\label{lem:p5smooth}
	Let $f$ be the solution of \eqref{eq:PludaDeTurck} constructed in \Cref{thm:STE}. If there exists $p>5$ and $\varepsilon>0$ such that $f \in W^{1,p}\left(\varepsilon,T_1; L^p(I;\R^d)\right)\cap L^p\left(\varepsilon,T_1; W^{4,p}(I;\R^d)\right)$ then $f \in \CalC^\infty((\varepsilon, T_1) \times I;\R^d)$.
\end{lem}
\begin{proof}
	See \cite[Chapter 3, Section 3.4]{Diss}.
\end{proof}

Now, \Cref{thm:Smoothing Main} is immediate.

\begin{proof}[{Proof of \Cref{thm:Smoothing Main}}]
	The solution $f$ in \Cref{thm:STE Main} is exactly the solution $f$ in \Cref{thm:STE}. By \Cref{prop:finMR20} we have $f\in W^{1,20}(\varepsilon,T_1;L^{20}(I;\R^d))\cap L^{20}(\varepsilon,T_1;L^{20}(I;\R^d))$ for any $0<\varepsilon<T_1$. Hence, by \Cref{lem:p5smooth}, we find that $f\in \CalC^{\infty}((\varepsilon,T_1)\times I;\R^d))$. 
\end{proof}

\section{Long time behavior and the {proof of \texorpdfstring{\Cref{thm:Convergence Main}}{Theorem 1.2}}}
\label{sec:Convergence}
In this section, we use the long time existence result in \cite{DPL17} to show the existence of a global solution of \eqref{eq:EF}. Moreover, we prove and use a \emph{refined {\L}ojasiewicz--Simon gradient inequality} to conclude convergence after reparametrization.

\subsection{Long-time existence after reparametrization}

As a first step towards proving \Cref{thm:Convergence Main}, we  establish long-time existence and subconvergence after reparametrization for our solution. The key ingredient is the smoothness of our solution and \cite[Theorem 1.1]{DPL17}.

\begin{thm}\label{thm:LTE}
	Let $f\in W^{1,2}(0,T;L^2(I;\R^d))\cap L^2(0,T;W^{4,2}(I;\R^d))$ be as in \Cref{thm:STE} and let $0<\varepsilon<T$. Then, there exists $\bar{\varepsilon}\in (\varepsilon, T)$ and $\hat{f}\in \CalC^{\infty}((0,\infty)\times I;\R^d)$ satisfying \eqref{eq:EF} such that
	\begin{enumerate}[(i)]
		\item $\hat{f}(t,x) = f(t,x)$ for all $0\leq t\leq \varepsilon, x\in I$;
		\item $\hat{f}(t, \cdot)$ has zero tangential velocity for all $t\geq \bar{\varepsilon}$;
		\item $\hat{f}$ subconverges smoothly as $t\to\infty$, after reparametrization with constant speed, to a constrained elastica, i.e.\ a solution  \eqref{eq:ElasticaFixedLength}.
	\end{enumerate}
\end{thm}

\begin{proof}
	By \Cref{thm:Smoothing Main}, the solution $f$ in \Cref{thm:STE} is instantaneously smooth. Thus, to simplify notation we may assume $f\in \CalC^{\infty}([\varepsilon, T]\times I;\R^d)$ for some $\varepsilon>0$ after possibly reducing $T>0$. Moreover, we may also assume a uniform bound from below on the arc-length element using \Cref{rem:STESolBounds}.
	
	Let $\theta\defeq \langle\partial_t {f}, \partial_{s_{{f}}} {f}\rangle$ be the tangential velocity of $f$. By the smoothness of $f$ and the bound on the arc-length element, the function $(t,r)\mapsto \frac{\theta(t,r)}{\abs{\partial_x f(t,r)}}$ is globally Lipschitz \blue{continuous} on $[\varepsilon, T]\times I$. For each $x\in I$, we consider the initial value problem
	\begin{align}\label{eq:ReparaODE}
	\left\lbrace\begin{array}{llr}
	\partial_t \Phi(t,x) &= - \frac{\theta(t,\Phi(t,x))}{\abs{\partial_x {f}(t, \Phi(t,x))}} \\
	\Phi (\varepsilon,x) &= x.
	\end{array}\right.
	\end{align}
	By classical ODE theory, there exists $\varepsilon< \hat{T}\leq T$ and a smooth family of reparame\-trizations $\Phi\colon [\varepsilon,\hat{T}]\times I\to I$ satisfying \eqref{eq:ReparaODE} and
	\begin{align}\label{eq:ODE Boundary}
	\Phi(t,y) &= y \quad \text{ for }t\in [\varepsilon,\hat{T}], y\in \partial I \nonumber\\
	\partial_x \Phi(t,x)&>0 \quad \text{ for all }(t,x)\in [\varepsilon,\hat{T}]\times I.
	\end{align} 
	Therefore, $\Phi(t, \cdot)$ is strictly increasing and a diffeomorphism of $I$ for each $t\in [\varepsilon, \hat{T}]$. A direct computation yields that the reparametrization ${f}_1(t,x)\defeq {f}(t, \Phi(t,x))$ satisfies
	\begin{align*}
	\partial_t {f}_1(t,x) &= \partial_t {f}(t,\Phi(t,x)) + \partial_x{f}(t, \Phi(t,x))\partial_t \Phi(t,x) \\
	&= \partial_t^{\perp}{f}(t,  \Phi(t,x)) + \theta(t,  \Phi(t,x)) \partial_{s_{{f}}}{f}(t, \Phi(t,x)) + \partial_x {f}(t, \Phi(t,x))\partial_t \Phi(t,x) \\
	&=\partial_t^{\perp}{f}(t,  \Phi(t,x)) \\
	&= -\nabla_{s_{f_1}}^2\vKap_{f_1}(t,x) -\frac{1}{2}\abs{\vKap_{f_1}(t,x)}^2\vKap_{f_1}(t,x) + \lambda({f_1})(t)\vKap_{f_1}(t,x),
	\end{align*}
	using that $f$ solves \eqref{eq:EF} and the transformation of the geometric quantities. For the boundary conditions, let $t\in [\varepsilon,\hat{T}]$, $y\in \partial I$ and note that $f_1(t,y)=f(t,y)=p_y$ and $ \partial_{s_{f_1}} f_1(t,y) =\partial_{s_f}f(t,y) =\tau_y$ by \eqref{eq:ODE Boundary}. Consequently, $f_1$ is a smooth solution of \eqref{eq:EF} on $[\varepsilon, \hat{T}]$ with tangential velocity zero and smooth initial datum $f(\varepsilon)$. By \cite[Theorem 1.1]{DPL17}, $f_1$ can be extended to a global smooth solution $\bar{f}$ on  $[\varepsilon, \infty)$ which subconverges, after reparametrization with constant speed, to a constrained elastica as $t\to\infty$. 
	
	In particular, we have the identity
	\begin{align}\label{eq:barf vs f Phi}
	\bar{f}(t,x) = f(t,\Phi(t,x))\quad \text{ for all }\varepsilon\leq t\leq \hat{T}.
	\end{align}
	Now, let $\varepsilon<\bar{\varepsilon}<\hat{T}$ and $\Psi\colon [0, \hat{T}]\times I \to I$ be a smooth family of reparametrizations with
	\begin{align}\label{eq:Psi Properties}
	\Psi(t,x) = x\quad\text{ for all }0\leq t\leq \varepsilon; \qquad 
	\Psi(t,x) =\Phi(t,x) \quad\text{ for all }\bar{\varepsilon}\leq t\leq \hat{T}.
	\end{align}
	The existence of such a $\Psi$ is proven in \Cref{lem:Repara Interpolation}. We now define
	\begin{align*}
	\hat{f}(t,x) \defeq \left\lbrace\begin{array}{ll}
	f(t,\Psi(t,x)) & \text{ for }0\leq t\leq \hat{T}, x\in I\\
	\bar{f}(t,x) & \text{ for } t\geq \bar{\varepsilon}, x\in I.
	\end{array}\right.
	\end{align*}
	Note that $\hat{f}$ is clearly smooth in $x$ for every $t\geq 0$ fixed. It is also smooth in $t$ for fixed $x\in I$, by \eqref{eq:barf vs f Phi} and \eqref{eq:Psi Properties}.
	Property (i) follows from \eqref{eq:Psi Properties}. Furthermore, by definition of $\hat{f}$ on $[\bar{\varepsilon},\infty)\times I$ we find that $\hat{f}=\bar{f}$ has zero tangential velocity and hence (ii) is satisfied. The last property follows since the asymptotic behavior of $\hat{f}$ is inherited from $\bar{f}$.
\end{proof}

\subsection{The length-preserving elastic flow as a gradient flow on a Hilbert manifold}
\label{subs:manifold}
{In this section, we show that the flow \eqref{eq:EF} is in fact a gradient flow on a suitable submanifold of curves.}

\begin{prop}Let $p_0,p_1 \in \R^d, \tau_0, \tau_1 \in \mathbb{S}^{d-1}$ and $\ell\in \R$ such that \eqref{eq:f_0assumption} holds. Then \begin{align*}
	\mathcal{X} \defeq \left\{ f\in W^{4,2}_{Imm}(I;\R^d)\mid f(y)=p_y \text{ and }\partial_s f(y) = \tau_y \text{ for } y\in \partial I, \CalL(f) = \ell\right\}.
	\end{align*} is a weak Riemannian splitting analytic submanifold of $W^{4,2}(I;\R^d)$ with codimension $4d-1$.
\end{prop}
\begin{proof}
	By the Sobolev embedding $W^{4,2}(I;\R^d)\hookrightarrow \mathcal{C}^{1}(I;\R^{d})$, the set of $W^{4,2}$-immersions denoted by $W^{4,2}_{Imm}(I;\R^d)$ is open in $W^{4,2}(I;\R^{d})$. The function
	\begin{align*}
		\mathcal{G}\colon W^{4,2}_{Imm}(I;\R^{d}) \to \R\times (\R^d)^2\times (\mathbb{S}^{d-1})^2, \mathcal{G}(f) \defeq \begin{pmatrix}
	\mathcal{L}(f) \\
	 f(0) \\ 
	 f(1) \\
	  \partial_{s_f}f(0)
	  \\ \partial_{s_f}f(1)
		\end{pmatrix}
	\end{align*}
	is an analytic map. Moreover, its differential is given by
	\begin{align*}
		&d \mathcal{G}_f \colon W^{4,2}(I;\R^d)\to \R\times (\R^d)^2 \times \mathcal{T}_{\partial_{s}f(0)}\mathbb{S}^{d-1}\times \mathcal{T}_{\partial_{s}f(1)}\mathbb{S}^{d-1},\\ 
		&d \mathcal{G}_f(u) = \begin{pmatrix}
		-\int_{I} \langle\vKap_{{f}}, u\rangle\diff s_{f} \\
		u(0) \\
		u(1) \\
		\frac{\partial_x u(0)}{\abs{\partial_x f(0)}} - \frac{\langle \partial_x u(0), \partial_x f(0)\rangle \partial_x f(0)}{\abs{\partial_x f(0)}^3} \\
		\frac{\partial_x u(1)}{\abs{\partial_x f(1)}} - \frac{\langle \partial_x u(1), \partial_x f(1)\rangle \partial_x f(1)}{\abs{\partial_x f(1)}^3}
		\end{pmatrix}
	\end{align*}
	for $f\in W^{4,2}_{Imm}(I;\R^{d})$ and $u \in W^{4,2}(I;\R^{d})$. It is not difficult to see, that $d\mathcal{G}_f$ is surjective if $f\in \mathcal{X} = \mathcal{G}^{-1}\left(\{(\ell, p_0, p_1, \tau_0, \tau_1)^T\}\right)$. Indeed, let $\alpha\in \R$, $q_y\in \R^{d}$, $z_y \in T_{\partial_{s}f(y)}\mathbb{S}^{d-1}$ for $y=0,1$. We have ${	\mathcal{T}_{\partial_s f(y)}\mathbb{S}^{d-1} = \{z\in \R^{d} \mid \langle z, \partial_x f(y)\rangle =0 \}}$. 
	Clearly, we can find an immersed curve $u\in W^{4,2}(I;\R^{d})$ with $u(y)=q_y$ and  $\frac{\partial_x u(y)}{\abs{\partial_x f(y)}} = z_y$ for $y=0,1$. Now, using the characterization of the tangent space, for $v\in \CalC^{\infty}_0(I;\R^d)$ we find
	\begin{align*}
		d\mathcal{G}_f(u+v) = \begin{pmatrix}
		-\int_{I} \langle\vKap_{{f}}, u+v\rangle\diff s_{f} \\
		q_0\\
		q_1\\
		z_0\\
		z_1\\
		\end{pmatrix},
	\end{align*}
	since adding $v$ does not change the boundary behavior.
	Moreover, as $\vKap_{f}\not \equiv 0$ using $f\in \mathcal{X}$ and \eqref{eq:f_0assumption}, we can choose $v$ such that $\int_{I}\langle \vKap_f, v\rangle\diff s_{f} =\varepsilon\neq 0$. Setting $\beta\defeq \int_{I}\langle\vKap_f, u\rangle\diff s_{f}$ and $w\defeq u -\frac{\alpha+\beta}{\delta} v$, we find $\int_{I} \langle \vKap_f, w\rangle\diff s_{f} = \beta - (\alpha + \beta) = -\alpha$, hence we have shown $d\mathcal{G}_f(w) = (\alpha, q_0, q_1, z_0, z_1)$, so $d\mathcal{G}_f$ is surjective.
	
	Consequently, $\mathcal{X}\subset W^{4,2}(I;\R^{d})$ is a splitting  submanifold by \cite[Theorem 3.5.4]{Abraham04} with codimension $1+2d+2(d-1)=4d-1$. Like in \cite{Fabian}, the analytic form of the implicit function theorem can be used to show that $\mathcal{X}$ is in fact analytic. The tangent space is given by 
	\begin{align}\label{eq:T_fX characterization}
		\mathcal{T}_f \mathcal{X} &= \ker d\mathcal{G}_f \nonumber \\
		&= \{u \in W^{4,2}(I;\R^{d})\mid u=0 \text{ on }\partial I, \partial_x^{\perp_f}u = 0\text{ on }\partial I, \int_{I}\langle \vKap_f, u\rangle\diff s_{f}=0\}.
	\end{align}
	Since \eqref{eq:EF} is a $L^2(\diff s_f)$ gradient flow, it is natural to endow $\mathcal{X}$ with the Riemannian metric $\langle u,v\rangle_{L^2(\diff s_f)} = \int_{I} \langle u,v\rangle \diff s_{f}$ for $u,v \in \mathcal{T}_f \mathcal{X}$. Note that since \blue{$\mathcal{T}_f\mathcal{X}$} is certainly not complete with respect to the induced norm, the metric is only \emph{weakly Riemannian} (cf.\ \cite[Definition 5.2.12]{Abraham04}).
\end{proof}

{
	It is not difficult to see that by \eqref{eq:T_fX characterization} the right hand side of the evolution \eqref{eq:EFEvolution} is the projection of the full $L^2
(\diff s_f)$-gradient $\nabla\CalE(f)$ onto the $L^2(\diff s_f)$-closure of the tangent space $\mathcal{T}_f\mathcal{X}$. This implies that \eqref{eq:EFEvolution} is the gradient flow of $\CalE$ on the manifold $\mathcal{X}$.
}

\subsection{The constrained \texorpdfstring{{\L}}{L}ojasiewicz--{S}imon gradient inequality}
In this subsection, we establish a {\L}ojasiewicz--Simon inequality for $\CalE$ on $\mathcal{X}$. To do so, we have to deal with the invariance of both energies \blue{$\mathcal{E}$ and $\mathcal{G}$}, which unfortunately creates large kernels for their first and second variations. Like in \cite{CFS09,Loja}, we work around this issue by restricting the energy to normal directions and using the implicit function theorem.

In the following, we always assume that the assumptions  \eqref{eq:BCCompatibility} and \eqref{eq:f_0assumption} are satisfied.

\begin{defi}
	Fix $\bar{f}\in \mathcal{X}$ and define $V_c\defeq W^{4,2}(I;\R^d)\cap  W^{2,2}_{0}(I;\R^d)$. 
	We define the \emph{space of normal vector fields along $\bar{f}$} by
	\begin{align*}
	W^{4,2,\perp}(I;\R^d) \defeq \{ f\in W^{4,2}(I;\R^d) \mid \langle f, \partial_x \bar{f}\rangle = 0 \text{ on } I\}.
	\end{align*}
	Moreover, we define  $H^{\perp}\defeq L^{2,\perp}(I;\R^d)\defeq \{ u\in L^{2}(I;\R^d)\mid \langle u, \partial_x \bar{f}
	\rangle
	= 0 \text{ a.e.}\}$ and $V_c^{\perp} \defeq V_c \cap W^{4,2,\perp}(I;\R^d)$. Both are Hilbert spaces and the $L^2$-orthogonal projection onto $H^{\perp}$ is given by the pointwise projection $P^{\perp}(f) \defeq f - \langle f, \partial_s \bar{f}\rangle\partial_s \bar{f}$.
	
	Moreover, by the embedding $W^{4,2}(I;\R^d)\hookrightarrow \mathcal{C}^{1}(I;\R^d)$ there exists $\varepsilon>0$ small enough such that for all $u\in W^{4,2,\perp}(I;\R^d)$ with $\norm{u}{W^{4,2}}<\varepsilon$, the curve $f=\bar{f}+u$ is immersed. Defining $U_\varepsilon \defeq \{ u\in V_c^{\perp} \mid \norm{u}{W^{4,2}}<\varepsilon\}$ we consider the energies
	\begin{align*}
	&L \colon U_\varepsilon\to \R, \quad L(u)=\CalL(\bar{f}+u) \text{ and }\\
	& E\colon U_\varepsilon\to \R, \quad E(u)=\CalE(\bar{f}+u).
	\end{align*}
\end{defi}
We have the following result.
\begin{prop}[{cf.\ \cite[Proof of Theorem 3.1, Remark 3.3]{Loja}}]\label{thm:ELojaAssum}
	The energy $E$ satisfies the following properties.
	\begin{enumerate}
		\item $E\colon U_\varepsilon\to \R$ is analytic;
		\item its gradient $\nabla E \colon U_\varepsilon \to H^{\perp}$ is analytic;
		\item the derivative $(\nabla E)^{\prime}(0)\colon V^{\perp}_{c} \to H^{\perp}$ is Fredholm with index zero.
	\end{enumerate}
\end{prop}
It is well known that this is sufficient to prove a {\L}ojasiewicz--Simon gradient inequality for $E$ (cf.\ \cite[Corollary 3.11]{Chill03}),\cite[Theorem 3.1]{Loja},\cite[Theorem 1.2]{Fabian}, \cite[Corollary 2.6]{Pozzetta}). However, in order to conclude a {\emph{constrained} or \emph{refined}} {\L}ojasiewicz--Simon gradient inequality, cf.\ (16) in \cite{Fabian}, we { also need to analyze the length functional.}
\begin{prop}\label{prop:GLojaAssum} The energy $L$ satisfies the following properties.
	\begin{enumerate}
		\item $L\colon U_\varepsilon\to \R$ is analytic.
		\item The gradient map $\nabla L \colon U_\varepsilon \to H^{\perp}$ is analytic.
		\item The derivative $(\nabla L)^{\prime}(0)\colon V_{c}^{\perp}\to H^{\perp}$ is compact.
		\item $L(0)=\ell$ and $\nabla L(0)\neq 0$.
	\end{enumerate}
\end{prop}
\begin{proof}
	\begin{enumerate}
		\item The map $U_{\varepsilon}\to \CalC(I;\R^d), u \mapsto \abs{\partial_x(\bar{f}+u)}$ is analytic by \cite[Lemma 3.4, 1.]{Loja}, and hence so is $L$.
		\item The $H^{\perp}$-gradient of $L$ is given by $\nabla L(u) = -P^{\perp} \left( \vKap_{\bar{f}+u}\abs{\partial_{x}(\bar{f}+u)}\right)$. Note that the map $U_{\varepsilon}\to L^{2}(I;\R^d), u \mapsto \vKap_{\bar{f}+u}$ is analytic by \cite[Lemma 3.4, 3.]{Loja}. Since the multiplication $L^2(I;\R^d)\times L^{\infty}(I;\R) \to L^2(I;\R^d), (f,\phi)\mapsto f\phi$ is analytic, so is the map $U_{\varepsilon}\to L^2(I;\R^d), u\mapsto  \vKap_{\bar{f}+u}\abs{\partial_{x}(\bar{f}+u)}$. The continuity and linearity of $P^{\perp}\colon L^2(I;\R^d)\to H^{\perp}$ yields the claim.
		\item We compute the second derivative using standard formulas for the variation of geometric quantities (see for instance \cite[Lemma 2.1]{DKS}). We have
		\begin{align*}
		(\nabla L)^{\prime}(0) u &= \dtzero \nabla L(tu) = -\dtzero P^{\perp} \left( \vKap_{\bar{f}+u}\abs{\partial_{x}(\bar{f}+u)}\right) \\
		&= -P^{\perp} \dtzero \vKap_{\bar{f}+tu} \abs{\partial_{x}\bar{f}} - P^{\perp} \vKap_{\bar{f}} \dtzero\abs{\partial_x (\bar{f}+tu)}\\
		&= -\left(\nabla_{s_{\bar{f}}}^2 u + \langle u, \vKap_{\bar{f}}\rangle \vKap_{\bar{f}}\right) \abs{\partial_{x}\bar{f}} + \vKap_{\bar{f}}\langle u, \vKap_{\bar{f}}\rangle \abs{\partial_x \bar{f}}.
		\end{align*}
		In particular, the operator $(\nabla L)^{\prime}(0)\colon V_c^{\perp}\to H^{\perp}$ is only of second order in $u$, hence compact by the Rellich--Kondrachov Theorem \cite[Theorem 7.26]{GilbargTrudinger}.
		\item $L(0)=\CalL(\bar{f})=\ell$ since $\bar{f}\in \mathcal{X}$. Since we have $\abs{\bar{f}(1)-\bar{f}(0)} = \abs{p_1-p_0} < \ell$, $\bar{f}$ cannot be part of a  straight line, hence $\vKap_{\bar{f}}\not\equiv 0$ and also $\abs{\partial_x \bar{f}}\neq 0$ since $\bar{f}$ is immersed.\qedhere
	\end{enumerate}
\end{proof}
This enables us to conclude the inequality in normal directions.
\begin{thm}\label{thm:LojaNormal}
	Suppose $\bar{f}\in \mathcal{X}$ is a constrained elastica. Then, there exist $C, \sigma>0$ and $\theta\in (0,\frac{1}{2}]$ such that for all $f=\bar{f}+u\in \mathcal{X}$ with $u\in V_c^{\perp}$ and $\norm{u}{W^{4,2}}\leq \sigma$ we have
	\begin{align*}
	\abs{\CalE(f)-\CalE(\bar{f})}^{1-\theta}\leq C\norm{\nabla_{L^2(\diff s_{f})} \CalE(f) + \lambda(f)\nabla_{L^2(\diff s_{f})} \CalL(f)}{L^2(\diff s_{f})}.
	\end{align*}
\end{thm}
\begin{proof}
	First, we verify the conditions of \cite[Corollary 5.2]{Fabian} for the energy $E$ and the constraint $\mathcal{G}(u)=L(u)-\ell$ on the spaces $V=V_{c}^{\perp}, H=H^{\perp}$. Note that $\nabla \mathcal{G} = \nabla{L}$. Clearly, $V_c^{\perp}\hookrightarrow H^{\perp}$ densely.
	Assumptions (ii) and (iii) follow from \Cref{thm:ELojaAssum}, whereas assumptions (iv)-(vi) are satisfied by \Cref{prop:GLojaAssum}. Note that $u=0$ is a constrained critical point of $E$ on $\mathcal{M}=\mathcal{G}^{-1}(\{0\})$ since $\bar{f}$ is a constrained elastica.
	
	Then, by \cite[Corollary 5.2]{Fabian} $E\vert_{\mathcal{M}}$ satisfies a constrained {\L}ojasiewicz--Simon gradient inequality, i.e.\ there exist $C, \sigma>0$ and $\theta\in (0,\frac{1}{2}]$ such that for all $u\in \mathcal{M}$ with $\norm{u}{W^{4,2}}\leq \sigma$ we have
	\begin{align*}
	\abs{E(u)-E(0)}^{1-\theta}\leq C \norm{P(u)\nabla E(u)}{L^2},
	\end{align*}
	where $P(u)\colon H^{\perp}\to H^{\perp}$ denotes the orthogonal projection onto the closure of the tangent space $\overline{T_u \mathcal{M}} = \{ y\in H^{\perp}\mid \langle \nabla L(u), y\rangle_{L^2} = 0\}$ (cf.\ \cite[Proposition 3.3]{Fabian}). Therefore, for \begin{align*}
	{\lambda}(f) = \frac{\langle\vKap_f, \nabla \CalE(f)\rangle_{L^2(\diff s_{f})}}{\norm{\vKap_f}{L^2(\diff s_{f})}^2}
	\end{align*}
	as in \eqref{eq:lambda} with $f=\bar{f}+u$ we have the estimate
	\begin{align*}
	\norm{P(u)\nabla E(u)}{L^2} = \norm{P(u)\left(\nabla E(u) + {\lambda} \nabla L(u)\right)}{L^2} \leq \norm{\nabla E(u) + {\lambda} \nabla L(u)}{L^2}.
	\end{align*}
	Moreover, we have $\nabla E(u) = \nabla_{L^2(\diff s_{f})} \CalE(f) \abs{\partial_x f}$ and $\nabla L(u) = \nabla_{L^2(\diff s_{f})} \CalL(f) \abs{\partial_x f}$. Consequently,
	\begin{align*}
	\norm{P(u)\nabla E(u)}{L^2} &\leq  \norm{P(u)\left(\nabla E(u) + {\lambda} \nabla L(u)\right)}{L^2} \\
	&\leq \norm{\nabla_{L^2(\diff s_{f})} \CalE(f) \abs{\partial_x f} + {\lambda} \nabla_{L^2(\diff s_{f})} \CalL(f) \abs{\partial_x f}}{L^2}\\
	&\leq \norm{\partial_x f}{L^{\infty}}^{\frac{1}{2}}\norm{\nabla_{L^2(\diff s_{f})} \CalE(f)  + {\lambda} \nabla_{L^2(\diff s_{f})} \CalL(f) }{L^2(\diff s_{f})}.
	\end{align*}
	Reducing $\sigma>0$ if necessary, we may assume that $\norm{\partial_x f}{L^{\infty}}$ is uniformly bounded for $\norm{f-\bar{f}}{W^{4,2}}\leq \sigma$ by the Sobolev embedding theorem. This proves the claim.	
\end{proof}
We use this to prove the full constrained {\L}ojasiewicz--Simon gradient inequality for not necessarily normal variations via the following reparametrization argument.
\begin{lem}[{\cite[Lemma 4.1]{Loja}}]\label{lem:Repara}
	Let $\bar{f}\in W^{5,2}(I;\R^d)$ be a regular curve. Then, there exists $\sigma>0$ such that for all $\psi\in V_c$ with $\norm{\psi}{W^{4,2}}\leq \sigma$, there exists a $W^{4,2}$-diffeomorphism $\Phi\colon I\to I$ such that 
	\begin{align}\label{eq:Repara}
	(\bar{f}+\psi)\circ \Phi = \bar{f}+\eta
	\end{align}
	for some $\eta\in V_c^{\perp}$. Moreover, given $\sigma>0$ there exists $\tilde{\sigma}=\tilde{\sigma}(\bar{f}, \sigma)>0$ such that for all $\psi\in V_c$ with $\norm{\psi}{W^{4,2}}\leq \tilde{\sigma}$ we have the above representation with $\norm{\eta}{W^{4,2}}\leq \sigma$.
\end{lem}
\begin{thm}\label{thm:LojaFull}
	Let $\bar{f}\in\mathcal{X}\cap W^{5,2}(I;\R^d)$ be a constrained elastica. Then there exist $C, \sigma>0$ and $\theta\in (0,\frac{1}{2}]$ such that 
	\begin{align*}
	\abs{\CalE({f})-\CalE(\bar{f})} \leq C  \norm{\nabla_{L^2(\diff s_{f})} \CalE(f) + \lambda(f)\nabla_{L^2(\diff s_{f})} \CalL(f)}{L^2(\diff s_{f})},
	\end{align*}
	for all $f\in \mathcal{X}$ such that $\norm{f-\bar{f}}{W^{4,2}}\leq \sigma$.
\end{thm}

\begin{proof}
	Let $C, \sigma>0, \theta\in (0,\frac{1}{2}]$ as in \Cref{thm:LojaNormal}, $\bar{f}\in \mathcal{X}$ be a constrained critical point of $\mathcal{E}$ on $\mathcal{X}$.
	By the regularity assumption on $\bar{f}$, we may use \Cref{lem:Repara}.

	 Thus, we find $\tilde{\sigma}>0$ such that \eqref{eq:Repara} holds for all $\psi\in V_c$ with $\norm{\psi}{W^{4,2}}\leq \tilde{\sigma}$ for some $\eta\in V_c^{\perp}$ with $\norm{\eta}{W^{4,2}}\leq \sigma$. Let $f\in \mathcal{X}$ such that $\norm{f-\bar{f}}{W^{4,2}}\leq \tilde{\sigma}$. Then by \Cref{lem:Repara}, there exist a diffeomorphism $\Phi\colon I \to I$ and $\eta\in V_c^{\perp}$ such that $f\circ \Phi= \bar{f}+\eta$. 
	
	Note that with $f, \bar{f}\in \mathcal{X}$ we also get $f\circ \Phi = \bar{f}+\eta\in \mathcal{X}$, since $\mathcal{L}(f) = \mathcal{L}(f\circ \Phi) = \ell$.	Since the elastic energy is invariant under reparametrization, we hence get using \Cref{thm:LojaNormal}
	\begin{align}\label{eq:LojaAllDirections}
	&\left\vert\mathcal{E}(f)-\mathcal{E}(\bar{f})\right\vert^{1-\theta}  =\left\vert\mathcal{E}(\bar{f}+\eta)-\mathcal{E}(\bar{f})\right\vert^{1-\theta} \nonumber\\
	&\qquad \leq C\norm{\nabla_{L^2(\diff s_{\bar{f}+\eta})} \CalE(\bar{f}+\eta) + \lambda(\bar{f}+\eta)\nabla_{L^2(\diff s_{\bar{f}+\eta})} \CalL(\bar{f}+\eta)}{L^2(\diff s_{\bar{f}+\eta})}.
	\end{align}
	Since $\lambda$ and the gradients are geometric, i.e transform correctly under reparametri\-zations, we have
	\begin{align*}
	\lambda(\bar{f}+\eta) &= \lambda(f\circ \Phi) = \lambda(f), \\
	\nabla_{L^2(\diff s_{\bar{f}+\eta})} \CalE(\bar{f}+\eta) &= \nabla_{L^2(\diff s_{f})} \CalE(f) \circ \Phi \\
	\nabla_{L^2(\diff s_{\bar{f}+\eta})} \CalL(\bar{f}+\eta) &= \nabla_{L^2(\diff s_{f})} \CalL(f) \circ \Phi.
	\end{align*}
	Consequently, we obtain
	\begin{align*}
	&\norm{\nabla_{L^2(\diff s_{\bar{f}+\eta})} \CalE(\bar{f}+\eta) + \lambda(\bar{f}+\eta)\nabla_{L^2(\diff s_{\bar{f}+\eta})} \CalL(\bar{f}+\eta)}{L^2(\diff s_{\bar{f}+\eta})} \\
	&\qquad = \norm{\nabla_{L^2(\diff s_{f})} \CalE(f)\circ \Phi + \lambda(f)\nabla_{L^2(\diff s_{f})} \CalL(f)\circ \Phi}{L^2(\diff s_{f\circ \Phi})} \\
	&\qquad =  \norm{\nabla_{L^2(\diff s_{f})} \CalE(f) + \lambda(f)\nabla_{L^2(\diff s_{f})} \CalL(f)}{L^2(\diff s_{f})}.
	\end{align*}
	Together with \eqref{eq:LojaAllDirections}, this implies the {\L}ojasiewicz--Simon gradient inequality for the elastic energy on $\mathcal{X}$.
\end{proof}

\subsection{Convergence}
In previous works (see e.g.\ \cite[p. 358 -- 359]{CFS09} and \cite[p. 2188 -- 2191]{Loja}), a lot of PDE theory and a priori parabolic Schauder estimates are needed to apply the {\L}ojasiewicz--Simon gradient inequality to conclude convergence for geometric problems. In this section, we introduce a novel inequality (see \Cref{lem:BeschtesLemma}) which enables us to significantly shorten this lengthy argument in the proof of \Cref{thm:Convergence Main}. We exploit the explicit structure of the constant speed reparametrization and the length bound to control the full velocity of the constant speed parametrization by the purely normal velocity of the original evolution. 
\begin{defi}\label{def:constSpeedRepara}
	Let $T\in (0, \infty]$ and let $f\colon[0,T)\times I\to\R^d$ be a family of immersed curves in $\R^d$.
	The \emph{constant speed $\CalL(f(t))$ reparametrization $\tilde{f}(t)$ of $f(t)$} is given by $\tilde{f}(t,x)\defeq f(t, \psi(t,x))$ where $\psi(t, \cdot)\colon{I}\to{I}$ is the inverse of $\varphi(t,\cdot)\colon{I}\to{I}$ given by
	\begin{align*}
	\varphi (t,x) := \frac{1}{\mathcal{L}(f(t))}\int_0^x |\partial_x f(t,z)|\diff z = \frac{1}{\mathcal{L}(f(t))}\int_0^x \diff s_{f(t)}.
	\end{align*}
\end{defi}

\begin{lem}
 \label{lem:BeschtesLemma}
	Suppose $T\in (0, \infty]$ and $f\colon[0,T)\times I\to\R^d$ is a family of curves in $\R^d$, such that $f(t,0)=p_0, f(t,1)=p_1$ and $\CalL(f(t)) >0$ for all $t\in (0,T]$. Then, if $\tilde{f}(t)$ is the constant speed $\CalL(f(t))$ reparametrization of $f(t)$,  for all $t\in [0,T)$ we have
	\begin{align*}
	\norm{\partial_t \tilde{f}(t)}{L^2(\diff x)}\leq \sqrt{\frac{2}{\mathcal{L}(f(t))}+16~\mathcal{E}(f(t))}\norm{\partial_t f}{L^2(\diff s_{f(t)})}.
	\end{align*}
{In particular, if $f$ evolves by the length preserving elastic flow \eqref{eq:EF}, we have \begin{align*}
	\norm{\partial_t \tilde{f}(t)}{L^2(\diff x)}\leq C\norm{\partial_t f}{L^2(\diff s_{f(t)})},
\end{align*}
for all $t\in (0,T]$, where $C=\sqrt{\frac{2}{\ell}+16\CalE(f_0)}$.}
\end{lem}
\begin{proof}
Recall that by \Cref{def:constSpeedRepara} we have 
\begin{align*}
\psi(t,\varphi (t,x)) = \varphi (t,\psi(t,x)) = x \text{ for all } t\in [0,T), x\in {I}.
\end{align*}
For the derivatives of $\varphi$ and $\psi$ we thus obtain 
\begin{enumerate}[(i)]
	\item $\partial_t \varphi (t,x)=-\frac{\partial_t \mathcal{L}(f(t))}{\mathcal{L}(f(t))^2} \int_0^{x} \diff s_{f(t)} - \frac{1}{\mathcal{L}(f(t))}\int_0^{x}\langle \partial_t f, \vec{\kappa}_{f(t)}\rangle \diff s_{f(t)}$;
	\item $\partial_x \varphi (t,x) = \frac{|\partial_x f(t,x)|}{\mathcal{L}(f(t))}$;
	\item[(iii)] $\partial_x \psi(t,\varphi (t,x)) = \left(\partial_x \varphi (t,x)\right)^{-1} = \frac{\mathcal{L}(f(t))}{|\partial_x f(t,x)|}$;
	\item[(iv)] $\begin{aligned}[t]
	&\partial_t \psi(t, \varphi (t,x)) = - \partial_x \psi(t,{\varphi(t,x)})~ \partial_t\varphi (t,x)\\
	&\qquad= \frac{\mathcal{L}(f(t))}{|\partial_x f(t,x)|}\left(\frac{\partial_t \mathcal{L}(f(t))}{\mathcal{L}(f(t))^2} \int_0^{x} \diff s_{f(t)} + \frac{1}{\mathcal{L}(f(t))}\int_0^{x}\langle \partial_t f, \vec{\kappa}_{f(t)}\rangle \diff s_{f(t)}\right).
	\end{aligned}$
\end{enumerate}
Now, we
estimate
\begin{align*}
\norm{\partial_t \tilde{f}(t)}{L^2(\diff x)}^2
&\leq 2\int_0^1 |(\partial_t f)(t, \psi(t,x))|^2\diff x \\
&\qquad + 2\int_0^1 |(\partial_x f)(t, \psi(t,x))|^2 ~|\partial_t \psi(t,x)|^2\diff x.
\end{align*}
Taking $y=\psi(t,x)$ and using $\psi(t,0)=0$, $\psi(t,1)=1$, we find
\begin{align*}
\norm{\partial_t \tilde{f}(t)}{L^2(\diff x)}^2 \leq& 2\Big(\int_{0}^{1}|\partial_t f(t,y)|^2 \frac{1} {\partial_x \psi(t, \varphi (t,y))}\diff y\\
&\quad+ |\partial_x f(t,y)|^2 ~ |\partial_t \psi (t, \varphi (t,y))|^2\frac{1} {\partial_x \psi(t, \varphi (t,y))}\diff y\Big)=: 2(A + B).
\end{align*}
For the first integral, we clearly have
\begin{align*}
A =  \int_0^1 |\partial_t f(t,y)|^2 \frac{|\partial_x f(t,y)|}{\mathcal{L}(f(t))}\diff y = \frac{1}{\mathcal{L}(f(t))}\norm{\partial_t f}{L^2(\diff s_{f(t)})}^2.
\end{align*}
For the second part, note that by (iv), we have
\begin{align*}
&B =\int_0^1 \Bigg\vert \frac{\partial_t \mathcal{L}(f(t))}{\mathcal{L}(f(t))} \int_0^{y} \diff s_{f(t)} + \int_0^{y}\langle \partial_t f, \vec{\kappa}_{f(t)}\rangle \diff s_{f(t)} \Bigg\vert^2 \frac{|\partial_x f(t,y)|}{\mathcal{L}(f(t))}\diff y.
\end{align*}
Now, using the boundary conditions, we have $\partial_t \CalL(f(t)) = -\int_{I}\langle\vKap_{f(t)}, \partial_t f(t)\rangle\diff s_{f(t)}$ and the Cauchy--Schwarz inequality yields
\begin{align*}
B &\leq 2 \int_0^1 \Bigg(\left\vert\frac{\partial_t \mathcal{L}(f(t))}{\mathcal{L}(f(t))} \int_0^y\diff s_{f(t)}\right\vert^2 + \left\vert \int_0^y \langle\partial_t f(t), \vec{\kappa}_{f(t)}\rangle \diff s_{f(t)}\right\vert^2 \Bigg)\frac{|\partial_x f(t,y)|}{\mathcal{L}(f(t))}\diff y \\
&\leq 2\int_0^1 \Bigg( \left(\int_0^1 |\langle \partial_t f(t), \vec{\kappa}_{f(t)}\rangle|\diff s_{f(t)}\right)^2 \\
&\qquad + \left(\int_0^1 |\langle \partial_t f(t), \vec{\kappa}_{f(t)}\rangle|\diff s_{f(t)}\right)^2\Bigg) \frac{|\partial_x f(t,y)|}{\mathcal{L}(f(t))}\diff y\\
&= 4 \left(\int_0^1 |\langle\partial_t f(t), \vec{\kappa}_{f(t)}\rangle |\diff s_{f(t)}\right)^2 \leq 4 \norm{\partial_t f(t)}{L^2(\diff s_{f(t)})}^2 \norm{\vec{\kappa}_{f(t)}}{L^2(\diff s_{f(t)})}^2 \\&= 8~ \mathcal{E}(f(t)) \norm{\partial_t f(t)}{L^2(\diff s_{f(t)})}^2. \qedhere
\end{align*}
\end{proof}

\begin{remark}
	Note that in the proof of \Cref{lem:BeschtesLemma}, we only used the boundary conditions to conclude that no boundary terms appear when integrating by parts. In particular, \Cref{lem:BeschtesLemma} also holds in the case of closed curves.
\end{remark}

Finally, we can prove our main convergence result.

\begin{proof}[{Proof of \Cref{thm:Convergence Main}}]
	Let $\varepsilon>0$ and let $\hat{f}\in \CalC^{\infty}((0, \infty)\times I;\R^d)$, $\bar{\varepsilon}>\varepsilon$ be as in \Cref{thm:LTE}. 
	The first statement of \Cref{thm:Convergence Main} follows from property (i) in \Cref{thm:LTE}, 
	and the fact that the solution $f$ in \Cref{thm:STE} lies in $\X_{T,2}\hookrightarrow\nolinebreak BUC([0,T];W^{2,2}(I;\R^d))$ by \Cref{prop:SobolevEmbeddings} and \eqref{eq:Besov=Sobolev}.
	
	For the convergence statement, let $\tilde{f}$ be the constant speed $\ell$ reparametrization of $\hat{f}$, cf.\ \Cref{def:constSpeedRepara}, and note that $\tilde{f}\in \CalC^{\infty}((0,\infty)\times I;\R^d)$. By \Cref{thm:LTE} (iii), there exists a sequence $t_n\to \infty$ and a smooth regular curve $f_{\infty}\colon I\to\R^{n}$, such that $\tilde{f}(t_n)\to f_\infty$ in $\mathcal{C}^{k}(I;\R^{n})$ for all $k\in \N_0$. Moreover, as a consequence of \Cref{thm:LTE}, $f_{\infty}$ is a smooth constrained elastica, i.e.\ a smooth solution of \eqref{eq:ElasticaFixedLength}.
	
	Recall from \Cref{thm:LTE} (ii) that $\hat{f}$ has tangential velocity zero for $t$ sufficiently large. Thus, we can without loss of generality assume $\mathcal{E}(\hat{f}(t))= \mathcal{E}(\tilde{f}(t)) >\mathcal{E}(f_\infty)$, since otherwise $\hat{f}(t)$ would be eventually constant by \eqref{eq:FlowMonotonical}, and hence convergent. 
	Moreover, since $\CalE(\hat{f}(t))$ is non increasing, we have that $\lim_{t\to\infty}\mathcal{E}(\hat{f}(t)) = \lim_{n\to\infty}\CalE(\tilde{f}(t_n))= \mathcal{E}(f_\infty)$.

	Since $f_{\infty}$ is smooth, by \Cref{thm:LojaFull}, there exists $\sigma, C_{LS}>0$ and $\theta\in (0,\frac{1}{2}]$ such that we have a refined {\L}ojasiewicz--Simon inequality, i.e.\ for all $g\in \mathcal{X}$ satisfying $\norm{g-f_{\infty}}{W^{4,2}}\leq \sigma$ we have
	\begin{align}\label{eq:LSfinfty}
	|\mathcal{E}(g)-\mathcal{E}(f_\infty)|^{1-\theta}\leq C_{LS} \norm{\nabla_{L^2(\diff s_{g})} \CalE(g) + \lambda(g)\nabla_{L^2(\diff s_{g})} \CalL(g)}{L^2(\diff s_{g})}.
	\end{align}
	Passing to a subsequence, we can assume $\norm{\tilde{f}(t_n,\cdot)-f_{\infty}}{W^{4,2}}<\sigma$ for all $n$. Define
	\begin{align*}
	s_n := \sup\left\lbrace s\geq t_n\mid \norm{\tilde{f}(t,\cdot)-f_{\infty}}{W^{4,2}}<\sigma \text{ for all } t\in[t_n, s]\right\rbrace
	\end{align*}
	and note that $s_n>t_n$ since $\tilde{f}$ is smooth.
	Define $G(t) := \left(\mathcal{E}(\tilde{f}(t))-\mathcal{E}(f_\infty)\right)^{\theta}$. By our assumption $\mathcal{E}(\tilde{f}(t))>\mathcal{E}(f_\infty)$, so we can compute on $[t_n, s_n)$  using that $\hat{f}$ solves \eqref{eq:EF} with $\theta\equiv 0$, so $\partial_t \hat{f} = -\nabla\CalE(\hat{f})-\lambda\nabla\CalL(\hat{f})$ and the fact that $\mathcal{E}$ is geometric, i.e.\ invariant under reparametrization
	\begin{align*}
	&-\frac{\diff}{\diff t}G  = \theta \left(\mathcal{E}(\tilde{f} ) - \mathcal{E}(f_\infty)\right)^{\theta-1}\left(-\frac{\diff}{\diff t}\mathcal{E}(\hat{f} )\right) \\
	&\quad= \theta \left(\mathcal{E}(\tilde{f} ) - \mathcal{E}(f_\infty)\right)^{\theta-1} \left(- \left\langle \nabla_{L^2(\diff s_{\hat{f} })}\mathcal{E}(\hat{f} ), \partial_t \hat{f} \right\rangle_{L^2(\diff s_{\hat{f} })}\right)\\
	&\quad = \theta\left(\mathcal{E}(\tilde{f} ) - \mathcal{E}(f_\infty)\right)^{\theta-1}  \norm{\nabla_{L^2(\diff s_{\hat{f} })} \CalE(\hat{f} ) + \lambda(\hat{f} )\nabla_{L^2(\diff s_{\hat{f} })} \CalL(\hat{f} )}{L^2(\diff s_{\hat{f} })} \norm{\partial_t \hat{f} }{L^2(\diff s_{\hat{f} })}.
	\end{align*}
	However, the quantity $\norm{\nabla_{L^2(\diff s_{\hat{f} })} \CalE(\hat{f} ) +\lambda(\hat{f} )\nabla_{L^2(\diff s_{\hat{f} })} \CalL(\hat{f} )}{L^2(\diff s_{\hat{f} })}$ is geometric, too. Thus
	\begin{align*}
	&-\frac{\diff}{\diff t}G  \\
	&= \theta\left(\mathcal{E}(\tilde{f} ) - \mathcal{E}(f_\infty)\right)^{\theta-1} \norm{\nabla_{L^2(\diff s_{\tilde{f} })} \CalE(\tilde{f} ) + \lambda(\tilde{f} )\nabla_{L^2(\diff s_{\tilde{f} })} \CalL(\tilde{f} )}{L^2(\diff s_{\tilde{f} })} \norm{\partial_t \hat{f} }{L^2(\diff s_{\hat{f} })}\\
	&
	\geq \frac{\theta}{C_{LS}} \norm{\partial_t \hat{f} }{L^2(\diff s_{\hat{f} })}.
	\end{align*}
	on $[t_n, s_n)$ by \eqref{eq:LSfinfty} and our choice of $s_n$. Therefore, by \Cref{lem:BeschtesLemma} we have
	\begin{align}\label{eq:dGdtUngl}
	-\frac{\diff}{\diff t}G(t) \geq {C} \norm{\partial_t \tilde{f}}{L^2(\diff x)},  
	\end{align}
	for all $t\in [t_n, s_n)$, where ${C} = C(\ell, \mathcal{E}(f_0), \theta, C_{LS})>0$.
	Let $t\in [t_n, s_n)$. Then
	\begin{align}\label{eq:goestozero}
	\norm{\tilde{f}(t)-\tilde{f}(t_n)}{L^2(\diff x)}\leq \int_{t_n}^{t}\norm{\partial_t \tilde{f}(\tau)}{L^2(\diff x)}\diff \tau \leq \frac{1}{{C}} G(t_n)\to 0
	\end{align}
	using \eqref{eq:dGdtUngl} and $\mathcal{E}(\tilde{f}(t_n)) \to\mathcal{E}(f_\infty)$ as $n\to \infty$. {We now assume that all of the $s_n$ are finite. Then, by continuity \eqref{eq:goestozero} also holds for $t=s_n$.
	 By the subconvergence result in \Cref{thm:LTE},
	passing to a subsequence we have $\tilde{f}(s_{n})\to \nolinebreak\psi$ smoothly as $n \to\nolinebreak \infty$. Moreover, by continuity and the definition of $s_n$, we have that $\norm{\psi-\nolinebreak f_\infty}{W^{4,2}}=\sigma$, whereas $\norm{\psi-f_\infty}{L^2(\diff x)} = \lim_{n\to\infty}\norm{\tilde{f}(s_{n})-\tilde{f}(t_{n})}{L^2(\diff x)} = 0$ by \eqref{eq:goestozero}, a contradiction.}
	
	Consequently, there has to exist some $n_0\in\N$ such that $s_{n_0}=\infty$, and this yields $\norm{\tilde{f}(t)-f_\infty}{W^{4,2}}<\sigma$ for all $t\geq t_{n_0}$. This means that \eqref{eq:dGdtUngl} holds for any $t\geq t_{n_0}$, thus $t\mapsto\norm{\partial_t \tilde{f}(t)}{L^2(\diff x)}\in L^1([0,\infty);\R)$. Hence, for all $t_{n_0}\leq t\leq t^{\prime}$ we have
	\begin{align*}
	\norm{\tilde{f}(t) - \tilde{f}(t^{\prime})}{L^2(\diff x)} \leq \int_{t}^{t^{\prime}}\norm{\partial_t\tilde{f}(\tau)}{L^2(\diff x)}\diff \tau \to 0,
	\end{align*}
	as $t,t'\to\infty$ by the dominated convergence theorem. Therefore,  $\lim_{t\to\infty}\tilde{f}(t)$ exists in $L^2(\diff x)$ and thus equals $f_\infty$. A subsequence argument shows that for any $k\in \N_0$ we have $\norm{\tilde{f}(t)-f_\infty}{\mathcal{C}^{k}(I;\R^d)}\to 0$ as $t\to\infty$, i.e.\ the convergence is smooth.
\end{proof}

\appendix
	\section{Explicit formulas in coordinates}
	In this section, we present the explicit representation of the geometric quantities appearing in this article. {They can be obtained by a straight forward calculation, see e.g.\ \cite[{(2.3)}]{GMP19}.}
	\begin{prop} Suppose $f\colon I \to\R^{d}$ is a smooth immersion. {With the arc-length element $\gamma = \abs{\partial_x f}$ } we have \label{prop:ZeugsInKoordinaten}
		\begin{enumerate}[(i)]
			\item $\vKap_f = \partial_s^{2} f  = \frac{\partial^2_x f}{\gamma^{2}} - \frac{\langle \partial_x^{2}f, \partial_x f\rangle}{\gamma^{4}} \partial_x f = \frac{(\partial_x^2 f)^{\perp}}{\gamma^4}$;
			\item $\nabla_{s_f} \vKap_f = \frac{\partial_x^{3}f}{\gamma^{3}} - \frac{\langle \partial_x^{3}f, \partial_x f\rangle}{\gamma^{5}}\partial_x f - 3 \frac{\langle\partial_x^{2} f, \partial_{x} f\rangle}{\gamma^{5}} \partial_x^{2} f + 3 \frac{\langle\partial_x^{2}f, \partial_{x}f\rangle^{2}}{\gamma^{7}}\partial_x f$;
			\item $\begin{aligned}[t]
				\nabla_{s_f}^2 \vKap_f &= \Big[ \frac{\partial_x^{4}f}{\gamma^{4}} -6 \frac{\langle\partial_{x}^{2} f, \partial_x f\rangle}{\gamma^{6}}\partial_x^{3}f - 4 \frac{\langle\partial_x^{3}f, \partial_x f\rangle}{\gamma^{6}}\partial_x^{2}f \\
				&\quad - 3 \frac{\abs{\partial_x^{2} f}^{2}}{\gamma^{6}} \partial_x^{2} f + 18 \frac{\langle\partial_x^{2}f, \partial_x f\rangle^{2}}{\gamma^{8}} \partial_x^{2}f\Big]^{\perp_f};
			\end{aligned}$
			\item $\begin{aligned}[t] \nabla\CalE(f) 
			&=\Big[\frac{\partial_x^{4}f}{\gamma^{4}} -6 \frac{\langle\partial_{x}^{2} f, \partial_x f\rangle}{\gamma^{6}}\partial_x^{3}f - 4 \frac{\langle\partial_x^{3}f, \partial_x f\rangle}{\gamma^{6}}\partial_x^{2}f \\
			&\quad - \frac{5}{2} \frac{\abs{\partial_x^{2} f}^{2}}{\gamma^{6}} \partial_x^{2} f + \frac{35}{2} \frac{\langle\partial_x^{2}f, \partial_x f\rangle^{2}}{\gamma^{8}} \partial_x^{2}f\Big]^{\perp_f}.\end{aligned}$
		\end{enumerate}
		Here $\nabla \CalE(f)$ denotes the $L^2(\diff s_f)$-gradient of $\CalE$ at $f$.
	\end{prop}

	\begin{lem}\label{lem:FormulasLambdaExplicit} Let $f$ be a smooth immersed curve with arc-length element $\gamma$. Then
		\begin{enumerate}[(i)]
			\item $\abs{\vKap_{f}}^{4} = \gamma^{-8} \abs{\partial_x^2 f}^{4} - 2\gamma^{-10} \abs{\partial_x^{2}f}^{2}\langle\partial_x^{2}f, \partial_x f\rangle^{2} + \gamma^{-12} \langle\partial_x^{2}f, \partial_x f\rangle^{4}$;
			\item $\begin{aligned}[t]
					\abs{\nabla_{s_f}\vKap_{f}}^{2}&= \gamma^{-6}\abs{\partial_x^{3} f}^{2} - \gamma^{-8} \langle\partial_x^{3}f, \partial_x f\rangle^{2} - 6\gamma^{-8}\langle\partial_x^{3}f, \partial_x^{2}f\rangle\langle\partial_x^{2}f, \partial_x f\rangle \\
				&\quad + 6 \gamma^{-10} \langle\partial_x^{3}f, \partial_x f\rangle \langle\partial_x^{2}f, \partial_x f\rangle^{2} + 9 \gamma^{-10}\langle\partial_x^{2}f, \partial_x f\rangle^{2} \abs{\partial_x^{2}f}^{2} \\
				&\quad - 9 \gamma^{-12}\langle\partial_x^{2}f, \partial_x f\rangle^{4};
			\end{aligned}$
			\item $\begin{aligned}[t]
			\langle \nabla_{s_f}\vKap_{{f}}, \vKap_f\rangle &= \gamma^{-5} \langle\partial_x^{3}f, \partial_x^{2}f\rangle - \gamma^{-7} \langle\partial_x^{3}f, \partial_x f\rangle\langle\partial_x^{2}f, \partial_x f\rangle \\
			&\quad - 3\gamma^{-7}\langle\partial_x^{2}f, \partial_x f\rangle \abs{\partial_x^{2} f}^2 + 3\gamma^{-9} \langle\partial_x^{2}f, \partial_x f\rangle^{3}.
			\end{aligned}$
		\end{enumerate}
	\end{lem}
	
	\section{Function spaces}
	
	In this section, we collect all relevant information on the function spaces for maximal $L^p$-regularity. Most of the embedding results are collected in \blue{\cite{DHP07}, see also} \cite{MS12}, where even a polynomial weight $t^{1-\mu}$ in time is allowed. As discussed in \Cref{rem:p=2 justify}, time weights do not allow to prove short-time existence with weaker initial data, and hence we restrict ourselves to the case $\mu=1$.
	
	Let $J\subset \R$ be an interval, $1\leq p<\infty$. For any $s\in (0,\infty)\setminus\N$  and a Banach space $E$, the ($E$-valued) \emph{Sobolev--Slobodetskii space} and the \emph{Bessel potential space}, respectively, are given by real and complex interpolation
	\begin{align*}
		W^{s,p}(J;E)&\defeq \left(W^{[s],p}(J;E), W^{[s]+1,p}(J;E)\right)_{s-[s],p};\\
		H^{s,p}(J;E)&\defeq \left(W^{[s],p}(J;E), W^{[s]+1,p}(J;E)\right)_{s-[s]};
	\end{align*}
	where $W^{k,p}(J;E)$ denotes the usual Bochner--Sobolev space for $k\in \N_0$. Recall from \cite[Theorem 2.4.1 (a), Definition 4.2.1]{TriebelInterpol} that the \emph{Besov spaces} are given by
	\begin{align*}
		B^{s}_{p,q}(I;\R^d) \defeq \left(W^{m_1,p}(I;\R^d), W^{m_2,p}(I;\R^d)\right)_{\theta,q},
	\end{align*}
	where $s=(1-\theta)m_1+\theta m_2$, $\theta\in (0,1), m_1, m_2\in\N_0, m_1<m_2, p,q\in [1,\infty)$. By \cite[Definition 2.3.1 (d) and Theorem 2.3.2 (d)]{TriebelInterpol} we have the relation 
	\begin{align}\label{eq:Besov=Sobolev}
	B^{s}_{p,p}(I;\R^d)&=W^{s,p}(I;\R^d) \quad \text{ for }s\not\in \N;\nonumber\\
	B^{s}_{2,2}(I;\R^d)&=W^{s,2}(I;\R^d) \quad \text{ for }s>0.
	\end{align}
	
	Moreover, recall from \cite[Section 2]{DHP07}, that in the setting of the maximal regularity spaces $\X_{T,p}$, the spaces of zeroth and first order boundary data are given by
	\begin{align}\label{eq:BoundarySpaces}
	\nonumber
	\CalD^0_{T,p} &\defeq W^{1-\frac{1}{4p},p}(0,T;L^p(\partial I;\R^d)) \cong   W^{1-\frac{1}{4p},p}(0,T;(\R^d)^2);  \\
	\CalD^1_{T,p} &\defeq W^{\frac{3}{4}-\frac{1}{4p},p}(0,T;L^p(\partial I;\R^d)) \cong W^{\frac{3}{4}-\frac{1}{4p},p}(0,T;(\R^d)^2).
	\end{align}
	
	\blue{We now recall the Sobolev embeddings for the maximal regularity space $\mathbb{X}_{T,p}$. We emphasize that the operator norms of the embeddings below might blow up as $T\to 0+$. However, in the solution space with vanishing trace, i.e.\ the space $ \prescript{}{0}{\X_{T,p}}$, they are bounded independently of $T$.}
		
	\begin{prop}\label{prop:SobolevEmbeddings}
		Let \blue{$0<T\leq\infty$, let $p\geq 2$, and let $k\in \{0,\dots,4\}$}.
		\begin{enumerate}[(i)]
			\item $\X_{T,p} \hookrightarrow BUC([0,T], B^{4(1-\frac{1}{p})}_{p,p})(I;\R^d))$ with the estimate
			\begin{align*}
			\norm{f}{BUC([0,T];B^{4(1-\frac{1}{p})}_{p,p}(I;\R^d))} \leq \blue{C(p)}\norm{f}{\X_{T,p}} \quad \text{for }f\in \prescript{}{0}{\X_{T,p}}.
			\end{align*}
			\item $\X_{T,p} \hookrightarrow \CalC^{\alpha}([0,T];\CalC^{1,\alpha}(I;\R^d))$ for some $\alpha\in (0,1)$ with the estimate
			\begin{align*}
			\norm{f}{\CalC^{\alpha}([0,T]; \CalC^{1,\alpha}(I;\R^d))}\leq \blue{C(p,\alpha)}  \norm{f}{\X_{T,p}}\quad \text{for }f\in \prescript{}{0}{\X_{T,p}}.
			\end{align*}
			\item The $k$-th spatial derivative is continuous as a map 
			\begin{align*}
			\partial_x^k \colon \X_{T,p} &\to H^{\frac{4-k}{4}, p}(0,T;L^p(I;\R^d)) \cap  L^p(0,T; H^{4-k,p}(I;\R^d)) \\
			&\hookrightarrow W^{\frac{(4-k)\theta}{4}, p}(0,T;W^{(4-k)(1-\theta),p}(I;\R^d)) \text{ for all } \theta\in (0,1),
			\end{align*}
			with the estimate
			\begin{align*}
			\norm{\partial_x^k f}{H^{\frac{4-k}{4}, p}(0,T;L^p(I;\R^d)) \cap  L^p(0,T; H^{4-k,p}(I;\R^d))}\leq \blue{C(k,p)} \norm{f}{\X_{T,p}}\quad \text{for }f\in \prescript{}{0}{\X_{T,p}}.
			\end{align*}
			\item The spatial trace of the $k$-th spatial derivative is continuous as a map
			\begin{align*}
			\tr_{\partial I}\partial_x^k \colon \X_T &\to W^{\frac{4-k}{4} - \frac{1}{8}, p}(0,T;L^p(\partial I;\R^d)) \cong W^{\frac{4-k}{4} - \frac{1}{8}, p}(0,T;(\R^d)^2),
			\end{align*}
			with the estimate 
			\begin{align*}
			\norm{ \tr_{\partial I}\partial_x^k f}{W^{\frac{4-k}{4} - \frac{1}{8}, p}(0,T;(\R^d)^2)}\leq \blue{C(k,p)} \norm{f}{\X_{T,p}}\quad \text{for }f\in \prescript{}{0}{\X_{T,p}}.
			\end{align*}
		\end{enumerate}
	
	\end{prop}
	\begin{proof}
		\blue{For $T=\infty$ and $I$ replaced by $\R$, the statements follow from the corresponding results in \cite[Section 3]{DHP07}. The statements in our case can then be obtained by considering appropriate temporal and spatial extension operators, see for instance \cite[Lemma 2.5 and (3.2)]{MS12}. When restricted to $\prescript{}{0}{\X_{T,p}}$, the operator norm of the temporal extension does not depend on $T$ by \cite[Lemma 2.5]{MS12}, which implies the above $T$-independent estimates.
		}
	\end{proof}
		
	\begin{remark}\label{rem:H=W}
		For a Hilbert space $E$ and $s\in (0,\infty),p=2$, the Bessel potential spaces coincide with the Slobodetskii spaces, i.e.\ $H^{s,2}(0,T;E) = W^{s,2}(0,T;E)$ with equivalence of norms, cf.\ \cite[Corollary 4.37]{Lunardi18}. A particular consequence of this is that in the case $p=2$ we get from \Cref{prop:SobolevEmbeddings} (iii) that \blue{for $k\in\N$, $k\leq 4$, and $T\in (0,1]$ we have}
		\begin{align*}
		\partial_x^k \colon \X_{T,2} &\to W^{\frac{4-k}{4}, 2}(0,T;L^2(I;\R^d)) \cap  L^2(0,T; W^{4-k,2}(I;\R^d))
		\end{align*}
		is continuous, with the estimate
		\begin{align*}
		\norm{\partial_x^k f}{W^{\frac{4-k}{4}, 2}(0,T;L^2(I;\R^d)) \cap  L^2(0,T; W^{4-k,2}(I;\R^d))}\leq \blue{C(k)} \norm{f}{\X_{T,2}}\quad \text{for }f\in \prescript{}{0}{\X_{T,2}}.
		\end{align*}
	\end{remark}
	A crucial tool in proving the contraction estimates in \Cref{subsec:Contraction} is the precise control of the integrability of the spatial derivatives and their spatial trace, with operator norm bounded independent of $T$. As in \Cref{subsec:Contraction}, we restrict to the case $p=2$ here.
	\begin{prop}\label{prop:SobolevEmbeddings2}
		Let \blue{$T\in (0,1]$, $k\in \N$, $k\leq 4$,} $\rho_1, \rho_2 \in [1, \infty)$.
		\begin{enumerate}[(i)]
			\item If there exists $\theta\in [0,1]$ such that $\frac{4-k}{4}\theta - \frac{1}{2}\geq -\frac{1}{\rho_1}$ and $(4-k)(1-\theta)-\frac{1}{2}\geq -\frac{1}{\rho_2}$ then $\partial_x^k \colon {\X_{T,2}}\to L^{\rho_1}(0,T;L^{\rho_2}(I;\R^d))$ with the estimate
			\begin{align*}
				\norm{\partial_x^k f}{L^{\rho_1}(0,T;L^{\rho_2}(I;\R^d))} \leq \blue{C(k,\theta, \rho_1, \rho_2)}\norm{f}{\X_{T,2}} \quad\text{ for all }f\in \prescript{}{0}{\X_{T,2}}.
			\end{align*}
			\item If $\frac{4-k}{4}-\frac{5}{8}\geq -\frac{1}{\rho_1}$, then $\tr_{\partial I} \partial_x^k \colon {\X_{T,2}}\to L^{\rho_1}(0,T;(\R^d)^2)$ with the estimate
			\begin{align*}
				\norm{\tr_{\partial I}\partial_x^k f}{L^{\rho_1}(0,T;(\R^d)^2} \leq \blue{C(k, \rho_1)}\norm{f}{\X_{T,2}} \quad\text{ for all }f\in \prescript{}{0}{\X_{T,2}}.
			\end{align*}
		\end{enumerate}
	\end{prop}
	\begin{proof}
		We first prove the estimates.
		\begin{enumerate}[(i)]
			\item Using first \Cref{prop:SobolevEmbeddings} (iii) and \Cref{rem:H=W}, then interpolation, and in the last line the usual Sobolev embedding both in the temporal and spatial variable, we find
			\begin{align}\label{eq:dx^k T_0 Sobolev embedding}
				\partial_x^k \colon \prescript{}{0}{\X_{\blue{1},2}} &\to W^{\frac{4-k}{4},2}(0,\blue{1};L^2(I;\R^d)) \cap L^2(0,\blue{1}, W^{4-k, 2}(I;\R^d)) \nonumber\\
				&\hookrightarrow W^{\frac{4-k}{4}\theta,2}(0,\blue{1};W^{(4-k)(1-\theta),2}(I;\R^d)) \nonumber\\
				& \hookrightarrow L^{\rho_1}(0, \blue{1}; L^{\rho_2}(I;\R^d)).
			\end{align}	
			Now, by \cite[Lemma 2.5]{MS12}, there exists an extension operator $E_T$ from $(0,T)$ to $(0,\blue{1})$ such that
			$E_T \colon \prescript{}{0}{\X_{T,2}} \to \prescript{}{0}{\X_{\blue{1},2}}$ has operator norm \blue{independent of $T$}. Then, for any $f\in \prescript{}{0}{\X_{T,2}}$ we have using \eqref{eq:dx^k T_0 Sobolev embedding}
			\begin{align*}
				\norm{\partial_x^k f}{L^{\rho_1}(0,T;L^{\rho_2}(I;\R^d))} &\leq  		\norm{\partial_x^k(E_T f)}{L^{\rho_1}(0,\blue{1};L^{\rho_2}(I;\R^d))} \\
				& \leq \blue{C(k, \theta, \rho_1, \rho_2)} \norm{E_T f}{\X_{\blue{1},2}} \\
				&\leq \blue{C(k, \theta, \rho_1, \rho_2)} \norm{f}{\X_{T,2}}.
			\end{align*}
			\item Since $\tr_{\partial I}\partial_k f$ only depends on the temporal variable, we first use \Cref{prop:SobolevEmbeddings} (iv) and \Cref{rem:H=W} and then the Sobolev embedding to find
			\begin{align}\label{eq:tr dx^k T_0 Sobolev embedding}
			\tr_{\partial I}\partial_x^k \colon \prescript{}{0}{\X_{\blue{1},2}} &\to W^{\frac{4-k}{4}-\frac{1}{8},2}(0,\blue{1};(\R^d)^2)\nonumber\\
			& \hookrightarrow L^{\rho_1}(0, \blue{1}; (\R^d)^2).
			\end{align}	
			Again, using the extension operator, we find for any $f\in \prescript{}{0}{\X_{T,2}}$
			\begin{align*}
			\norm{\tr_{\partial I}\partial_x^k f}{L^{\rho_1}(0,T;(\R^d)^2)} &\leq  		\norm{\tr_{\partial I}\partial_x^k(E_T f)}{L^{\rho_1}(0,\blue{1};(\R^d)^2)} \\
			& \leq\blue{C(k, \rho_1)} \norm{E_T f}{\X_{\blue{1},2}} \\
			&\leq \blue{C(k, \rho_1)} \norm{f}{\X_{T,2}}.
			\end{align*}
		\end{enumerate}
		The mapping properties follow from \eqref{eq:dx^k T_0 Sobolev embedding} and \eqref{eq:tr dx^k T_0 Sobolev embedding}.
	\end{proof}

	\section{Details of the contraction estimates}\label{app:contr}
	First, the following definition describes the structure of the nonlinearities in \eqref{eq:PludaDeTurck} which guarantees the desired contraction properties.
	\begin{defi}\label{def:nonlinearities}
		Let $(a,b)\in \N_0^{2}$. We denote by $A^{(a,b)}$ the set of bounded multilinear maps
		\begin{align}\label{eq:multilin varphi}
			\varphi\colon (\R^d)^{m}\times (\R^d)^{a}\times (\R^d)^{b}\to \R^{w}
		\end{align}
		for some $w\in \N$, $m\in \N_0$. Then, we define the set $\mathcal{A}^{(a,b)}$ of {\emph{multilinear maps of type $(a,b)$}} as the set of all maps $f\mapsto\Phi(f)$ acting via
		\begin{align*}
			&\Phi(f)(t,x) \\
			& = \varphi\Big(\underbrace{\partial_x f(t,x), \dots, \partial_x f(t,x)}_{m\text{-times}}, \underbrace{\partial_x^{2}f(t,x), \dots, \partial_x^{2}f(t,x)}_{a \text{-times}}, \underbrace{\partial_x^{3}f(t,x), \dots, \partial_x^{3}f(t,x)}_{b\text{-times}}\Big),
		\end{align*}
		for almost every $(t,x)\in (0,T)\times I$ where $\varphi\in A^{(a,b)}$.
	\end{defi}
	
	{
		\begin{remark}
			Note that we do not keep track of $m$, the number of first order derivatives appearing in $\Phi\in \CalA^{(a,b)}$. This is justified since by \Cref{prop:SobolevEmbeddings} (ii), the derivatives of first order of $f\in \X_{T}$ are in $\CalC([0,T]\times I;\R^d)$, and hence do not affect the integrability of $\Phi(f)$.
		\end{remark}
	}
	
	\begin{example}
		The map $f\mapsto\Phi(f) = {\langle\partial_x^2f, \partial_x f\rangle} \partial_x^3 f$ is in $\mathcal{A}^{(1,1)}$, since the derivatives of second and third order only appear linearly.
	\end{example}
	
	The following proposition yields for which parameters $(a,b)$ we get a contraction. Note that nonlinearities with \blue{this} structure appear in $\tilde{F}$ in \eqref{eq:defA} and $\lambda$ in \eqref{eq:lambdaIBP}. \blue{As in \Cref{subsec:Contraction}, we assume $T,M\leq 1$ and set $\mathbb{X}_T = \mathbb{X}_{T,2}$ to simplify notation.}
	
	\begin{prop}\label{thm:FundamentalContractionEstimates}
		Let $q\in (0,1)$ and let $\Phi\in \CalA^{(a,b)}$. Then, for \blue{$T=T(q,\bar{f}),M=M(q,\bar{f})\in (0,1] $} small enough, each of the following nonlinear maps is a well-defined $q$-contraction, i.e.\ Lipschitz \blue{continuous} with Lipschitz constant $q$.
		\begin{enumerate}[(i)]
			\item $\blue{\bar{B}_{T,M}}\to L^2(0,T;L^2), f\mapsto\Phi(f)$, if  $(a,b)=(1,1)$ or $(a,b)=(3,0)$.
			\item $\blue{\bar{B}_{T,M}} \to L^2(0,T), f\mapsto \int_{I} \Phi(f)\diff x$, if  $(a,b)=(0,2), (a,b)=(2,1)$ or $(a,b)=(4,0)$. 
			\item $\blue{\bar{B}_{T,M}} \to L^2(0,T;(\R^{d})^2), f\mapsto \tr_{\partial I} \Phi(f)$, if $(a,b)=(1,1)$ or $(a,b)=(3,0)$.
			\item $\blue{\bar{B}_{T,M}} \to L^2(0,T;L^2), f\mapsto\left(\gamma_0^{-4}-\gamma^{-4}\right)\partial_x^4f$.
		\end{enumerate}
	\end{prop}
	The following general functional analytic result gives sufficient conditions for a multilinear map to be a $q$-contraction for $T>0$ small. It \blue{is} the key ingredient in the proof of \Cref{thm:FundamentalContractionEstimates}.

	\begin{lem}\label{lem:MultilinBoundedLip}
		Let $1\leq q_1\leq \infty$ and suppose $(f_1, \dots, f_r)\mapsto\mu(f_1, \dots, f_r)$ is a multilinear map such that for all $f_1, \dots, f_r\in \X_{T,2}$ we have
		\begin{align}\label{eq:MultiLinBounded}
			\norm{\mu(f_1, \dots, f_r)}{L^{q_1}(0,T;Z)}\leq C \prod_{j=1}^{r} \norm{S\partial_x^{d_j} f_j}{X_j}.
		\end{align}
		Here, we have $d_1, \dots, d_r \in \{0,\dots, 3\}$, $S\in \{\Id, \tr_{\partial I}\}$ and $Z,X_1, \dots, X_r$ are Banach spaces such that \blue{there exists $C\in (0,\infty)$ independent of $T$ with}
		\begin{enumerate}[(i)]
			\item {$\partial_x^{d_i} \colon \X_T \to X_i$ and for $f\in \prescript{}{0}{\X_T}$ we have $\norm{S \partial_x^{d_i}f}{X_i}\leq \blue{C}\norm{f}{\X_T}$ for all $i=1, \dots,r$.}
			\item for all $j=1, \dots r$ one of the following conditions is satisfied.		\begin{enumerate}
				\item There exists $\alpha>0$ with $\norm{S\partial_x^{d_j}f}{X_j}\leq \blue{C} T^{\alpha} \norm{f}{\X_T}$ for all $f\in\prescript{}{0}{\X}_T$.
				\item There exists $k\neq j$ with $\norm{S\partial_x^{d_k}f}{X_k}\to 0$ as $T\to 0$ for all $f\in \X_T$.
			\end{enumerate}
		\end{enumerate}
		Then, setting $\mu(f) = \mu(f, \dots, f)$, we have $\mu(f)\in L^{q_1}(0,T;Z)$ for all $f\in \X_T$ and for any $q\in (0,1)$, there exist \blue{$M = M(q,r,\bar{f}),T = T(q,r, \bar{f})\in (0,1]$} small enough, such that for all $f,\tilde{f}\in \blue{\bar{B}_{T,M}}$ we have
		\begin{align*}
			\norm{\mu(f)-\mu(\tilde{f})}{L^{q_1}(0,T;Z)}\leq q\norm{f-\tilde{f}}{\X_T}.
		\end{align*}
	\end{lem}
	{
		\begin{remark}\label{rem:p<infty L2.8 ii b}
			When applying \Cref{lem:MultilinBoundedLip}, we always work with Banach spaces of the type $X_j=L^{p_j}(0,T;L^{q_j})$ and $Z=L^{p_0}$, for some $p_0, p_j, q_j \in [1,\infty]$. Note that (ii) b) is always satisfied if there exist $k\neq j$ with $p_k<\infty$, since then $\lim_{T\to 0} \norm{f}{L^{p_k}(0,T;L^{q_k})}\to 0$ by dominated convergence.
		\end{remark}
	}
	\begin{proof}[{Proof of \Cref{lem:MultilinBoundedLip}}]
		Let $f, \tilde{f}\in \blue{\bar{B}_{T,M}}$. Adding and subtracting zeroes and using the multilinearity, we get
		\begin{align*}
			\mu(f)-\mu(\tilde{f}) &= \mu(f-\tilde{f},f,\dots, f) + \mu(\tilde{f}, f-\tilde{f}, f, \dots, f)\\
			&\quad + \dots + \mu(\tilde{f}, \dots, f-\tilde{f}, f)+\mu(\tilde{f}, \dots, \tilde{f}, f-\tilde{f}).
		\end{align*}
		Thus, using \eqref{eq:MultiLinBounded}, we get
		\begin{align}
			\label{eq:MultiLinBounded2}
			&\norm{\mu(f)-\mu(\tilde{f})}{L^{q_1}(0,T;Z)}\nonumber \\
			&\leq C\sum_{j=1}^r   \norm{S\partial_x^{d_1}\tilde{f}}{X_1} \dots 
			\norm{S \partial_x^{d_j}(f-\tilde{f})}{X_j}
			\dots\norm{S\partial_x^{d_r}f}{X_r}.
		\end{align}
		We now show that the contraction property is valid for each summand in \eqref{eq:MultiLinBounded2}. Note that for all $k\in \{1, \dots,r\}$ by (i) we have
		\begin{align}\label{eq:NormfBounded0}
			\norm{S\partial_x^{d_k}f}{X_k} &\leq \norm{S\partial_x^{d_k}(f-\bar{f})}{X_k} + \norm{S\partial_x^{d_k}\bar{f}}{X_k}\nonumber \\
			&\leq \blue{C}\norm{f-\bar{f}}{\X_T} + \norm{S\partial_x^{d_k}\bar{f}}{X_k} \leq \blue{C}\left( M + \norm{S\partial_x^{d_k}\bar{f}}{X_k}\right).
		\end{align}
		In particular, for \blue{$T\leq 1$, $M\leq 1$} we find
		\begin{align}\label{eq:fNormBounded}
			\norm{S\partial_x^{d_k}f}{X_k}, \norm{S\partial_x^{d_k}\tilde{f}}{X_k}\leq \blue{C(\bar{f})}.
		\end{align}
		
		Now, let $j\in \{1, \dots, r\}$. If (ii) a) is satisfied, using $f(0)=\tilde{f}(0)=f_0$, we find
		\begin{align*}
			&\norm{S\partial_x^{d_1}\tilde{f}}{X_1} \dots \norm{S\partial_x^{d_{j-1}}\tilde{f}}{X_{j-1}}\norm{ S\partial_x^{d_j}(f-\tilde{f})}{X_j}\norm{S\partial_x^{d_{j+1}}f}{X_{j+1}}\dots\norm{S\partial_x^{d_r}f}{X_r}\\ &
			\quad \leq \blue{C(\bar{f})} T^{\alpha}\norm{f-\tilde{f}}{\X_T} \leq \frac{q}{Cr} \norm{f-\tilde{f}}{\X_T},
		\end{align*}
		for \blue{$T = T(\alpha, \bar{f})>0$} small enough. Otherwise, if (ii) b) is satisfied, we estimate using (i) for the $j$-th factor,  \eqref{eq:NormfBounded0} for the $k$-th factor and \eqref{eq:fNormBounded} for the remaining factors, to get
		\begin{align*}
			&\norm{S\partial_x^{d_1}\tilde{f}}{X_1} \dots \norm{S\partial_x^{d_{j-1}}\tilde{f}}{X_{j-1}}\norm{S \partial_x^{d_j}(f-\tilde{f})}{X_j}\norm{S\partial_x^{d_{j+1}}f}{X_{j+1}}\dots\norm{S\partial_x^{d_r}f}{X_r}\\
			&\quad\leq \blue{C(\bar{f})} \left(M + \norm{S\partial_x^{d_k}\bar{f}}{X_k}\right) \norm{f-\tilde{f}}{\X_T}.
		\end{align*}
		By (ii) b), $\lim_{T\to 0}\norm{S\partial_x^{d_k}\bar{f}}{X_k} = 0$. Consequently, for $\blue{T = T(q,r,\bar{f}), M = M(q,r, \bar{f}) \in (0,1]}$ small enough we find
		\begin{align*}
			&\norm{S\partial_x^{d_1}\tilde{f}}{X_1} \dots \norm{S\partial_x^{d_{j-1}}\tilde{f}}{X_{j-1}}\norm{S \partial_x^{d_j}(f-\tilde{f})}{X_j}\norm{S\partial_x^{d_{j+1}}f}{X_{j+1}}\dots\norm{S\partial_x^{d_r}f}{X_r} \\
			&\qquad \leq \frac{q}{Cr}\norm{f-\tilde{f}}{\X_T}.
		\end{align*}
		All in all, we have proven
		\begin{align*}
			\norm{\mu(f)-\mu(\tilde{f})}{L^{q_1}(0,T;Z)}\leq q \norm{f-\tilde{f}}{\X_T}\quad \text{ for }f, \tilde{f}\in \blue{\bar{B}_{T,M}}.&\qedhere
		\end{align*}
	\end{proof}
	
%
	
	Together with the embedding results in \Cref{prop:SobolevEmbeddings} and \Cref{prop:SobolevEmbeddings2}, we can now prove \Cref{thm:FundamentalContractionEstimates}. 
	
	\begin{proof}[{Proof of \Cref{thm:FundamentalContractionEstimates}}]
		Let $f, \tilde{f}\in \blue{\bar{B}_{T,M}}\subset \X_T$ with $\blue{T=T(q,\bar{f}), M=M(q,\bar{f})\in (0,1]}$ small enough such that \Cref{lem:boundsgamma} is satisfied. { The strategy for the proof of cases (i)-(iii) is to apply  \Cref{lem:MultilinBoundedLip}. To that end, we use H\"older's inequality in time and space, and then verify the assumptions of \Cref{lem:MultilinBoundedLip} using \Cref{prop:SobolevEmbeddings2}.
			In the following, we denote by $f_1,\dots,f_m$,  $g,g_1, g_2, g_3, g_4, h, h_1, h_2$ general functions in $\blue{\bar{B}_{T,M}}$.}
		
		\underline{Case (i):} 	If  $(a,b)=(1,1)$, {by H\"older's inequality we have}
		\begin{align} \label{eq:Phi_i_11}
			&\norm{\varphi(\partial_x f_1, \dots, \partial_x f_m, \partial_x^2g, \partial_x^3h)}{L^2(0,T;L^2)}\nonumber\\
			&\quad \leq C(\varphi)\prod_{j=1}^m\norm{\partial_x f_j}{\infty}\norm{\partial_x^2g}{L^{4}(0,T;L^{8})} \norm{\partial_x^3 h}{L^{4}(0,T;L^{\frac{8}{3}})}.
		\end{align}
		Using \Cref{prop:SobolevEmbeddings} (ii) we have
		\begin{align*}
			&	\partial_x \colon \X_T \to \CalC^{\alpha}([0,T]; \CalC(I;\R^d)),\nonumber \\
			&	\text{with the estimate } \norm{\partial_x f}{\CalC^{\alpha}([0,T];\CalC(I;\R^d))}\leq \blue{C}\norm{f}{\X_T} \text{ for }f\in \prescript{}{0}{\X_{T}}.
		\end{align*}	
		Therefore, 
		we find 
		\begin{align}\label{eq:dxC_alpha}
			&	\partial_x \colon \X_T \to \CalC([0,T]\times I;\R^d),\nonumber \\
			&	\text{with the estimate } \norm{\partial_x f}{\infty}\leq \blue{C}T^{\alpha}\norm{f}{\X_T}  \text{ for }f\in \prescript{}{0}{\X_{T}},
		\end{align}		
		such that (ii) a) in \Cref{lem:MultilinBoundedLip} is satisfied. 
		Next, using \Cref{prop:SobolevEmbeddings2} (i) with $k=2$ and $\theta=\frac{3}{4}$ yields 
		\begin{align}\label{eq:dx^2L8L4}
			&	\partial_x^2 \colon \X_T \to L^8(0,T;L^4), \nonumber\\
			&	\text{with the estimate } \norm{\partial_x^2 f}{L^8(0,T;L^4)}\leq \blue{C}\norm{f}{\X_T} \quad\text{ for all }f\in \prescript{}{0}{\X_{T}},
		\end{align}
		since $\frac{4-2}{4}\cdot\frac{3}{4}-\frac{1}{2}\geq -\frac{1}{8}$ and $(4-2)(1-\frac{3}{4}) - \frac{1}{2}\geq -\frac{1}{4}$. Similarly for the third derivative with $\theta=\frac{1}{2}$ we get
		\begin{align}\label{eq:dx^3L8/3L4}
			&	\partial_x^3 \colon \X_T \to L^{\frac{8}{3}}(0,T;L^4), \nonumber\\
			&	\text{with the estimate } \norm{\partial_x^3 f}{L^{\frac{8}{3}}(0,T;L^4)}\leq \blue{C} \norm{f}{\X_T} \quad\text{ for all }f\in \prescript{}{0}{\X_{T}}.
		\end{align}
		Thus, condition (ii) a) in \Cref{lem:MultilinBoundedLip} is satisfied for the fist $m$ factors in \eqref{eq:Phi_i_11} by \eqref{eq:dxC_alpha}, whereas for the remaining factors condition (ii) b) holds. More precisely, for $j=m+1$ choosing $k=m+2$ works and conversely $j=m+2,k=m+1$, using \Cref{rem:p<infty L2.8 ii b}. 
		
		{The case $(a,b)=(3,0)$ can be treated similarly, using H\"older to obtain} 
		\begin{align*}
			&\norm{\varphi(\partial_x f_1, \mydots, \partial_x f_m, \partial_x^2 g_1, \partial_x^2 g_2, \partial_x^2 g_3)}{L^2(0,T;L^2)} \leq C\prod_{j=1}^m\norm{\partial_x f_j}{\infty} \prod_{j=1}^{3}\norm{\partial_x^2g_j}{L^{6}(0,T;L^{6})},
		\end{align*}
		and then \Cref{prop:SobolevEmbeddings2} (i) with $k=2$, $\theta=\frac{2}{3}$ to get
		\begin{align}\label{eq:dx^2L6L6}
			&	\partial_x^2 \colon \X_T \to L^6(0,T;L^6), \nonumber\\
			&	\text{with the estimate } \norm{\partial_x^2 f}{L^6(0,T;L^6)}\leq \blue{C} \norm{f}{\X_T}\quad\text{ for all }f\in \prescript{}{0}{\X_{T}}.
		\end{align}
		
		\underline{Case (ii):} First, we have the following basic estimate
		\begin{align*}
			\Norm{\int_I \Phi(f)\diff x-\int_{I}\Phi(\tilde{f})\diff x}{L^2(0,T)} \leq \norm{\Phi(f)-\Phi(\tilde{f})}{L^2(0,T;L^1)}.
		\end{align*}
		{It hence suffices to show that $\X_T\to L^2(0,T;L^1), f\mapsto\Phi(f)$ is a $q$-contraction. To that end, we  use \Cref{lem:MultilinBoundedLip} with $Z=L^1, S=\Id$.} 
		
		If $(a,b)=(0,2)$, we have by H\"older's inequality
		\begin{align}\label{eq:Phi_ii_02}
			&\norm{\varphi(\partial_xf_1, \dots \partial_x f_m, \partial_x^3h_1, \partial_x^3h_2)}{L^2(0,T;L^1)} \leq C \prod_{j=1}^m\norm{\partial_x f_j}{\infty} \prod_{j=1}^{2}\norm{\partial_x^{3}h_j}{L^{4}(0,T;L^2)}.
		\end{align}
		Now, using \Cref{prop:SobolevEmbeddings2} (i) with $k=3$ and $\theta =1$, we have
		\begin{align}\label{eq:dx^3L4L2}
			&	\partial_x^3 \colon \X_T \to L^{4}(0,T;L^2(I;\R^d))\nonumber\\
			&	\text{with the estimate } \norm{\partial_x^3 f}{L^4(0,T;L^2)}\leq \blue{C} \norm{f}{\X_T} \quad\text{ for all }f\in \prescript{}{0}{\X_{T}}.
		\end{align}
		Consequently, the last two factors in \eqref{eq:Phi_ii_02} satisfy condition (i) and (ii) b) in \Cref{lem:MultilinBoundedLip}, cf.\ \Cref{rem:p<infty L2.8 ii b}. For the first $m$ factors, we may once again use \eqref{eq:dxC_alpha} to deduce that conditions (i) and (ii) a) in \Cref{lem:MultilinBoundedLip} are satisfied. 
		
		If $(a,b)=(2,1)$, we proceed similarly, first using H\"older to get
		\begin{align*}
			&\norm{\varphi(\partial_x f_1, \dots \partial_xf_m,\partial_x^2g_1, \partial_x^2g_2, \partial_x^3h)}{L^2(0,T;L^1)} \\
			&\quad \leq C \prod_{j=1}^m\norm{\partial_x f_j}{\infty} \prod_{j=1}^{2}\norm{\partial_x^2g_j}{L^{8}(0,T;L^4)}\norm{\partial_x^3 h}{L^4(0,T;L^2)},
		\end{align*}
		and then applying \eqref{eq:dxC_alpha}, \eqref{eq:dx^2L8L4} and \eqref{eq:dx^3L4L2}.
		For $(a,b)=(4,0)$, we may apply H\"older's inequality to obtain
		\begin{align*}
			&\norm{\varphi(\partial_xf_1, \dots \partial_xf_m,\partial_x^2g_1, \partial_x^2g_2, \partial_x^2g_3,\partial_x^2 g_4)}{L^2(0,T;L^1)}  \\
			&\quad \leq C\prod_{j=1}^m\norm{\partial_x f_j}{\infty} \prod_{j=1}^{4}\norm{\partial_x^2g_j}{L^{8}(0,T;L^4)},
		\end{align*}
		and then use \eqref{eq:dxC_alpha} and \eqref{eq:dx^2L8L4}.
		
		\underline{Case (iii):} Again, we use \Cref{lem:MultilinBoundedLip}, now with $Z=(\R^d)^2$ and $S=\tr_{\partial I}$. If $(a,b)=(1,1)$ we obtain by H\"older's inequality
		\begin{align}\label{eq:Phi_iii_11}
			&\norm{\tr_{\partial I} \varphi(\partial_xf_1, \dots, \partial_xf_m, \partial_x^2g,\partial_x^3h)}{L^2(0,T;(\R^d)^{2})} \nonumber\\
			&\quad \leq C\prod_{j=1}^m\norm{\tr_{\partial I}\partial_x f_j}{\infty} \norm{\tr_{\partial I}\partial_x^2g}{L^8(0,T;(\R^d)^2)}\norm{\tr_{\partial I}\partial_x^3 h}{L^{\frac{8}{3}}(0,T;(\R^d)^2)}.
		\end{align}
		Note that by \Cref{prop:SobolevEmbeddings2} (ii), we have
		\begin{align}\label{eq:trdx^2L8}
			&	\tr_{\partial I}\partial_x^2 \colon \X_T \to  L^{8}(0,T;(\R^d)^2), \nonumber\\
			&	\text{with the estimate } \norm{\tr_{\partial I}\partial_x^2 f}{L^{8}(0,T;(\R^d)^2)}\leq \blue{C} \norm{f}{\X_T}\quad\text{ for all }f\in \prescript{}{0}{\X_{T}},
		\end{align}
		whereas for the third derivative, we obtain
		\begin{align}\label{eq:trdx^3L8/3}
			&	\tr_{\partial I}\partial_x^3 \colon \X_T \to L^{\frac{8}{3}}(0,T;(\R^d)^2),\nonumber \\
			&	\text{with the estimate } \norm{\tr_{\partial I}\partial_x^3 f}{L^{\frac{8}{3}}(0,T;(\R^d)^2)}\leq \blue{C} \norm{f}{\X_T}\quad\text{ for all }f\in \prescript{}{0}{\X_{T}}.
		\end{align}
		{As in cases (i) and (ii), we then use the mapping properties and the estimates in \eqref{eq:dxC_alpha}, \eqref{eq:trdx^2L8} and \eqref{eq:trdx^3L8/3} together with \Cref{rem:p<infty L2.8 ii b} to verify that the assumptions of \Cref{lem:MultilinBoundedLip} are satisfied.}
		
		If $(a,b)=(3,0)$,  we proceed similarly, first using H\"older to obtain
		\begin{align*}
			&\norm{\tr_{\partial I} \varphi(\partial_xf_1, \dots, \partial_xf_m, \partial_x^2g_1, \partial_x^2g_2, \partial_x^2g_3)}{L^2(0,T;(\R^d)^{2})} \nonumber\\
			&\quad \leq C\prod_{j=1}^m\norm{\partial_x f_j}{\infty} \prod_{j=1}^{3}\norm{\partial_x^2g_j}{L^6(0,T;(\R^d)^2)},
		\end{align*}
		and then \eqref{eq:dxC_alpha} for the first order terms and \Cref{prop:SobolevEmbeddings2} (ii) with $k=2$, yielding
		\begin{align*}
			&\tr_{\partial I}\partial_x^2 \colon \X_T\to L^{6}(0,T;(\R^d)^2), \\
			&\text{with the estimate } \norm{\tr_{\partial I}\partial_x^2 f}{L^{6}(0,T;(\R^d)^2)}\leq \blue{C} \norm{f}{\X_T}\quad\text{ for all }f\in \prescript{}{0}{\X_{T}}.
		\end{align*}
		
		\underline{Case (iv):} Let $q\in (0,1)$. For $f, \tilde{f}\in \blue{\bar{B}_{T,M}}$, we have
		\begin{align*}
			\abs{(\gamma_0^{-4}-\gamma^{-4})\partial_x^4 f - (\gamma_0^{-4}-\tilde{\gamma}^{-4})\partial_x^4 \tilde{f}} &\leq \abs{\gamma_0^{-4}-\gamma^{-4}} \abs{\partial_x^4 f-\partial_x^4 \tilde{f}} + \abs{\tilde{\gamma}^{-4}-\gamma^{-4}}\abs{\partial_x^4 \tilde{f}}.
		\end{align*}
		Thus, we may estimate
		\begin{align}\label{eq:Phi_iv}
			&\norm{(\gamma_0^{-4}-\gamma^{-4})\partial_x^4 f - (\gamma_0^{-4}-\tilde{\gamma}^{-4})\partial_x^4 \tilde{f}}{L^2(0,T;L^2)}\nonumber \\
			&\quad \leq \norm{\gamma_0^{-4}-\gamma^{-4}}{\infty}\norm{\partial_x^4 f-\partial_x^4 \tilde{f}}{L^2(0,T;L^2)} + \norm{\tilde{\gamma}^{-4}-\gamma^{-4}}{\infty} \norm{\partial_x^4 \tilde{f}}{L^2(0,T;L^2)}.
		\end{align}
		For the first term, we use the mean value theorem and \Cref{lem:boundsgamma} to conclude
		\begin{align}\label{eq:gamma_MWS}
			\abs{\gamma_0(x)^{-4}-\gamma(t,x)^{-4}} 
			\leq 4 \left(\frac{\inf_{I}\gamma_0}{2} \right)^{-5}\abs{\partial_x f_0(x) -\partial_x f(t,x)}.
		\end{align}
		Consequently, using \Cref{prop:SobolevEmbeddings} (ii) we may estimate
		\begin{align*}
			\norm{\gamma_0^{-4}-\gamma^{-4}}{\infty}
			&\leq C(\gamma_0) \sup_{(t,x)\in [0,T]\times I} \abs{\partial_x f(0,x)-\partial_xf(t,x)} \\
			&\leq C(\gamma_0) \sup_{(t,x)\in [0,T]\times I} \norm{f(0)-f(t)}{\CalC^{1+\alpha}(I;\R^d)} \\
			&\leq C(\gamma_0) \sup_{(t,x)\in [0,T]\times I} t^{\alpha} \norm{f}{\CalC^{\alpha
				}([0,T]; \CalC^{1+\alpha}(I;\R^d))} \\
			&\leq C(\gamma_0) T^{\alpha}\left( \norm{f-\bar{f}}{\CalC^{\alpha
				}([0,T]; \CalC^{1+\alpha}(I;\R^d))} +\norm{\bar{f}}{\CalC^{\alpha
				}([0,T]; \CalC^{1+\alpha}(I;\R^d))}\right) \\
			&\leq \blue{C(\gamma_0)} T^{\alpha} \left(\norm{f-\bar{f}}{\X_T}+\norm{\bar{f}}{\CalC^{\alpha
				}([0,T]; \CalC^{1+\alpha}(I;\R^d))}\right) \\
			&\leq \blue{C(\bar{f})} T^{\alpha}.
		\end{align*}
		Combined with the simple estimate $\norm{\partial_x^4 f -\partial_x^{4}\tilde{f}}{L^2(0,T;L^2)}\leq \norm{f-\tilde{f}}{\X_T}$ this yields a $\frac{q}{2}$-contraction estimate for the first part of \eqref{eq:Phi_iv}, taking $\blue{T=T(q,
		\bar{f})\in(0,1]}$ small enough. For the remaining part, we use \eqref{eq:gamma_MWS} with $\gamma_0$ replaced by $\tilde{\gamma}$ to conclude
		\begin{align*}
			\norm{\tilde{\gamma}^{-4}-\gamma^{-4}}{\infty}\leq 4 \left(\inf_{I}\frac{\gamma_0}{2}\right)^{-5} \norm{\partial_x f-\partial_x \tilde{f}}{\infty} \leq \blue{C(\bar{f})} T^{\alpha}\norm{f-\tilde{f}}{\X_T},
		\end{align*}
		and $\norm{\partial_x^4 f}{L^2(0,T;L^2)} \leq  \norm{f-\bar{f}}{\X_T} +  \norm{\partial_x^4 \bar{f}}{L^2(0,T;L^2)} \leq \blue{C(\bar{f})}$. Consequently, if $\blue{T=T(q, \bar{f})\in(0,1]}$ is small enough, the second part of \eqref{eq:Phi_iv} is a $\frac{q}{2}$-contraction.
	\end{proof}
	
	It is not difficult to see that the statement of \Cref{thm:FundamentalContractionEstimates} remains true if one allows multiplication by powers of the arc-length element.
	
	\begin{cor}\label{cor:FundamentalContracion}
		Let $q\in (0,1),\ell\in \N$, {$\Phi\in \CalA^{(a,b)}$}. For $\blue{T=T(q,\ell) \in (0,1]}$, $\blue{M=M(q, \ell)\in (0,1]}$ small enough, each of the following maps is a well-defined $q$-contraction.
		\begin{enumerate}[(i)]
			\item $\blue{\bar{B}_{T,M}}\to L^2(0,T;L^2), f\mapsto \gamma^{-\ell}\Phi(f)$, if $(a,b)=(1,1)$ or $(a,b)=(3,0)$.
			\item $\blue{\bar{B}_{T,M}} \to L^2(0,T), f\mapsto \int_{I} \gamma^{-\ell}\Phi(f)\diff x$, if  $(a,b)=(0,2), (a,b)=(2,1)$ or $(a,b)=(4,0)$. 
			\item $\blue{\bar{B}_{T,M}} \to L^2(0,T(\R^{d})^2), f\mapsto \tr_{\partial I}\gamma^{-\ell} \Phi(f)$, if  $(a,b)=(1,1)$ or $(a,b)=(3,0)$.
		\end{enumerate}
	\end{cor}
	\begin{proof}
		\underline{Well-definedness:} By \Cref{lem:boundsgamma} we can estimate $\abs{\gamma^{-\ell} \Phi(f)}\leq \frac{\inf\gamma_0}{2} \abs{\Phi(f)}$ for all $T,M>0$ small enough. Thus $f\mapsto\gamma^{-\ell}\Phi(f)$ maps into the correct space by \Cref{lem:MultilinBoundedLip}.
		
		\underline{Contraction:} Let $q\in (0,1)$ and let $f, \tilde{f}\in \blue{\bar{B}_{T,M}}$. For the first case, taking $T,M>0$ small enough, we have
		\begin{align*}
			&\norm{\gamma^{-\ell}\Phi(f)-\tilde{\gamma}^{-\ell}\Phi(\tilde{f})}{L^2(0,T;L^2)} \\
			&\quad \leq \norm{\gamma^{-\ell}-\tilde{\gamma}^{-\ell}}{\infty} \norm{\Phi(f)}{L^2(0,T;L^2)} + \norm{\tilde{\gamma}^{-\ell}}{\infty} \norm{\Phi(f)-\Phi(\tilde{f})}{L^2(0,T;L^2)} \\
			&\quad \leq \norm{\gamma^{-\ell}-\tilde{\gamma}^{-\ell}}{\infty} \left(\norm{\Phi(f)- \Phi(\bar{f})}{L^2(0,T;L^2)}+\norm{\Phi(\bar{f})}{L^2(0,T;L^2)}\right) \\
			&\quad \quad + \left( \norm{\tilde{\gamma}^{-\ell}-\bar{\gamma}^{-\ell}}{\infty}+\norm{\bar{\gamma}^{-\ell}}{\infty}\right) \norm{\Phi(f)-\Phi(\tilde{f})}{L^2(0,T;L^2)} \\
			&\quad \leq \blue{C(\bar{f})}\norm{\gamma^{-\ell}-\tilde{\gamma}^{-\ell}}{\infty}+ \left( \norm{\tilde{\gamma}^{-\ell}-\bar{\gamma}^{-\ell}}{\infty}+\norm{\bar{\gamma}^{-\ell}}{\infty}\right) q_2\norm{f-\tilde{f}}{\X_T} 
		\end{align*}
		using \Cref{thm:FundamentalContractionEstimates} for $q_2\in (0,1)$ to be chosen. With similar estimates one finds
		\begin{align*}
			&\Norm{\int_I \gamma^{-\ell}\Phi(f)\diff x - \int_I\tilde{\gamma}^{-\ell}\Phi(\tilde{f})\diff x}{L^2(0,T)} \\
			&\quad  \leq\blue{C(\bar{f})}\norm{\gamma^{-\ell}-\tilde{\gamma}^{-\ell}}{\infty}+ \left( \norm{\tilde{\gamma}^{-\ell}-\bar{\gamma}^{-\ell}}{\infty}+\norm{\bar{\gamma}^{-\ell}}{\infty}\right) q_2\norm{f-\tilde{f}}{\X_T}  
		\end{align*}
		and
		\begin{align*}
			&\norm{\tr_{\partial I}\gamma^{-\ell}\Phi(f) -\tr_{\partial I}\tilde{\gamma}^{-\ell}\Phi(\tilde{f})}{L^2(0,T;(\R^d)^2)} \\
			&\quad  \leq \blue{C(\bar{f})}\norm{\gamma^{-\ell}-\tilde{\gamma}^{-\ell}}{\infty}+ \left( \norm{\tilde{\gamma}^{-\ell}-\bar{\gamma}^{-\ell}}{\infty}+\norm{\bar{\gamma}^{-\ell}}{\infty}\right) q_2\norm{f-\tilde{f}}{\X_T}.
		\end{align*}
		
		We now prove that for any $q\in (0,1)$ the map  $\blue{\bar{B}_{T,M}}\to L^{\infty}((0,T)\times I)$, $f\mapsto \gamma^{-\ell}$ is a $q$-contraction for $T,M>0$ small enough. We find as in \eqref{eq:gamma_MWS}
		\begin{align*}
			\norm{\gamma^{-\ell}-\tilde{\gamma}^{-\ell}}{\infty} 
			&\leq \blue{C(\ell, \bar{f})} \norm{f-\tilde{f}}{\CalC^0([0,T];\CalC^{1}(I;\R^d))}\\
			&\leq \blue{C(\ell, \bar{f})} T^{\alpha}\norm{f-\tilde{f}}{\X_T}\leq \frac{q}{2} \norm{f-\tilde{f}}{\X_T},
		\end{align*}
		for $\blue{T=T(q, \ell, \bar{f})\in(0,1]}$ small enough using \Cref{prop:SobolevEmbeddings} (ii) and the fact that $f(0)=\tilde{f}(0)=f_0$. Thus, we find
		\begin{align*}
			&\blue{C(\bar{f})}\norm{\gamma^{-\ell}-\tilde{\gamma}^{-\ell}}{\infty}+ \left( \norm{\tilde{\gamma}^{-\ell}-\bar{\gamma}^{-\ell}}{\infty}+\norm{\bar{\gamma}^{-\ell}}{\infty}\right) q_2\norm{f-\tilde{f}}{\X_T} \\
			&\quad \leq \frac{q}{2}\norm{f-\tilde{f}}{\X_T} + \left(M + \norm{\bar{\gamma}^{-\ell}}{\infty}\right)q_2 \norm{f-\tilde{f}}{\X_T} \leq q\norm{f-\tilde{f}}{\X_T},
		\end{align*}
		choosing first $\blue{q_2 = q_2 (q, \ell, \bar{f}) \in (0,1)}$ sufficiently small and passing to a smaller $\blue{T = T(q, \ell, \bar{f})}$ and  $\blue{M= M(q, \ell, \bar{f})\in(0,1]}$ if necessary.
	\end{proof}
	
	\begin{proof}[Proof of \Cref{lem:NContraction}]
		First, taking $\blue{T=T(\bar{f}),M\in (0,1]}$ small enough such that \Cref{lem:boundsgamma,lem:lowerBoundsEnergy} hold, all terms are defined almost everywhere.  We observe that $\CalF(f)$ is a sum of terms as in \Cref{cor:FundamentalContracion} (i) and \Cref{thm:FundamentalContractionEstimates} (iv) by \eqref{eq:defA}, hence well-defined and a $q$-contraction for all $q\in (0,1)$, if $\blue{T=T(q, \bar{f})}$, $\blue{M=M(q, \bar{f})\in(0,1]}$ are small enough. 
		
		For $\Lambda$ we need to do one {additional} estimate. {For $f\in \blue{\bar{B}_{T,M}}$ and $T,M>0$ the scalar valued function $\lambda$ is in $L^{2}(0,T)$, since by \Cref{lem:lowerBoundsEnergy} the energy $\CalE(f)$ in the denominator of $\lambda$ (cf.\ \Cref{subsec:lambda}) is bounded from below uniformly in $t$}, whereas the nominator $N(f)$ is in $L^2(0,T)$ by \Cref{cor:FundamentalContracion} (ii) and (iii) and by the explicit formulas in \Cref{lem:FormulasLambdaExplicit} and \eqref{eq:lambdaIBP}. The term $\vKap_{{f}}$ is in $L^{\infty}(0,T)$ by the embedding $\X_T\hookrightarrow BUC([0,T]; W^{2,2}(I;\R^d))$, cf.\ \Cref{prop:SobolevEmbeddings} (i),\eqref{eq:Besov=Sobolev} and \Cref{prop:ZeugsInKoordinaten}.
		
		{Now, the crucial step is the proof of the contraction estimate for $\Lambda$. To that end,} let $f, \tilde{f}\in \blue{\bar{B}_{T,M}}$. Then, writing $\lambda(f)=\frac{N(f)}{2\CalE(f)}$ {as in \Cref{subsec:lambda}}, we find for almost every $(t,x)$
		\begin{align}\label{eq:ContractionLambda0}
			&\abs{\lambda(f)(t)\vKap_f(t,x) - \lambda(\tilde{f})(t)\vKap_{{\tilde{f}}}(t,x)}\nonumber\\
			&\quad \leq \abs{\lambda(f)(t)-\lambda(\tilde{f})(t)} \abs{\vKap_{{f}}(t,x)} + \abs{\lambda(\tilde{f})(t)} \abs{\vKap_f(t,x)-\vKap_{\tilde{f}}(t,x)} \nonumber\\
			&\quad \leq \frac{1}{2\CalE(f(t))\CalE(\tilde{f}(t))} \abs{N(f)(t)} \abs{\CalE(f(t))-\CalE(\tilde{f}(t))}\abs{\vKap_{{f}}(t,x)}\nonumber \\
			&\quad \quad + \frac{1}{2\CalE(f(t))}\abs{N(f)(t)-N(\tilde{f})(t)}\abs{\vKap_f(t,x)} + \abs{\lambda(f)(t)} \abs{\vKap_{{f}}(t,x)-\vKap_{{\tilde{f}}}(t,x)} \nonumber\\
			&\quad \leq C(f_0) \abs{N(f)(t)} \abs{\CalE(f(t))-\CalE(\tilde{f}(t))}\abs{\vKap_{{f}}(t,x)} \nonumber\\
			&\quad \quad + C(f_0)\abs{N(f)(t)-N(\tilde{f})(t)}\abs{\vKap_f(t,x)} + \abs{\lambda(f)(t)} \abs{\vKap_{{f}}(t,x)-\vKap_{{\tilde{f}}}(t,x)},
		\end{align}
		using that by \Cref{lem:lowerBoundsEnergy} the elastic energy is bounded from below. Taking the $L^2L^2$-norm in \eqref{eq:ContractionLambda0}, we are left with three terms. The first one is
		\begin{align}\label{eq:ContractionLambdaSplit1}
			&\norm{ \abs{N(f)} \abs{\CalE(f)-\CalE(\tilde{f})}\abs{\vKap_{{f}}}}{L^2(0,T;L^2)}\nonumber \\
			&\quad \leq \norm{N(f)}{L^2(0,T)} \norm{\CalE(f)-\CalE(\tilde{f})}{L^{\infty}(0,T)} \norm{\vKap_f}{L^{\infty}(0,T;L^2)}.
		\end{align}
		{Now, note that $N(f)$ is a sum of terms as in \Cref{cor:FundamentalContracion} (ii) and (iii) by \eqref{eq:lambdaIBP} and the explicit formulas in \Cref{lem:FormulasLambdaExplicit}. 
			Therefore, for any $q\in (0,1)$, we have
			\begin{align}\label{eq:N_f contraction}
				\norm{N(f) - N(\tilde{f})}{L^2(0,T)} \leq q \norm{f-\tilde{f}}{\X_{T}},
			\end{align}	if we take $\blue{T=T(q, \bar{f}),M=M(q, \bar{f})\in(0,1]}$ small enough. In particular, we can assume that $f\mapsto N(f)$ is $1$-Lipschitz.}
		
		For the elastic energy term, note that $\CalE$ is analytic, hence $\mathcal{C}^{1}$ on the space of $W^{2,2}$-immersions, cf.\ \Cref{thm:ELojaAssum}, in particular it is locally Lipschitz \blue{continuous} in a neighborhood of $f_0\in W^{2,2}_{Imm}(I;\R^d)$. Hence, there exists $C(f_0)>0$ such that $\abs{\CalE(h)-\CalE(\tilde{h})}\leq C(f_0)\norm{h-\tilde{h}}{W^{2,2}(I;\R^d)}$ for all $h$ and $\tilde{h}$ satisfying $\norm{h-f_0}{W^{2,2}}<\delta$ and $\norm{\tilde{h}-f_0}{W^{2,2}}\leq \delta$. 
		
		By \Cref{prop:SobolevEmbeddings}(i), we have $\prescript{}{0}{\X}_T \hookrightarrow BUC(0,T;W^{2,2})$ \blue{with operator norm independent of $T\in (0,1]$}. Consequently, we have
		\begin{align*}
			\norm{f(t)-f_0}{W^{2,2}}&\leq \norm{f(t)-\bar{f}(t)}{W^{2,2}}+ \norm{\bar{f}(t)-f_0}{W^{2,2}} \\
			&\leq \blue{C}M +\norm{\bar{f}(t)-\bar{f}(0)}{W^{2,2}} \leq \delta
		\end{align*}
		for $\blue{ T=T(\bar{f}), M\in (0,1]}$ small enough, and similarly $\norm{\tilde{f}(t)-f_0}{W^{2,2}}\leq \delta$. But then, using \Cref{prop:SobolevEmbeddings}(i), we have the estimate
		\begin{align}\label{eq:ELipschitz}
			\norm{\CalE(f)-\CalE(\tilde{f})}{L^{\infty}(0,T)}\leq C(f_0) \norm{f-\tilde{f}}{L^{\infty}(0,T;W^{2,2})} \leq \blue{C(\bar{f})} \norm{f-\tilde{f}}{\X_T}.
		\end{align}
		
		For the curvature term $\vKap_f$, note that $W^{2,2}_{Imm}(I;\R^d)\to L^2(I;\R^d), f\mapsto \vKap_{f} = \partial_{s_f}^2 f$ is analytic (cf.\ \Cref{thm:ELojaAssum}), in particular Lipschitz \blue{continuous} near $f_0$. The same argument as above yields
		\begin{align}\label{eq:KappaLipschitz}
			\norm{\vKap_{f}-\vKap_{\tilde{f}}}{L^{\infty}(0,T;L^2)}\leq C(f_0) \norm{f-\tilde{f}}{L^{\infty}(0,T;W^{2,2})} \leq \blue{C(\bar{f})}\norm{f-\tilde{f}}{\X_T}.
		\end{align}
		Now, we estimate
		\begin{align}\label{eq:KappaBounded}
			\norm{\vKap_{f}}{L^{\infty}(0,T;L^{2})}&\leq \norm{\vKap_{f}-\vKap_{\bar{f}}}{L^{\infty}(0,T;L^2)} + \norm{\vKap_{\bar{f}}}{L^{\infty}(0;T;L^2)}\nonumber \\
			& \leq \blue{C(\bar{f})} \norm{f-\bar{f}}{\X_T}+\norm{\vKap_{\bar{f}}}{L^{\infty}(0;T;L^2)}\leq \blue{C(\bar{f})}
		\end{align}
		and using \eqref{eq:N_f contraction}, we obtain the bound
		\begin{align}\label{eq:N_fBounded}
			\norm{N(f)}{L^2(0,T)}&\leq \norm{N(f)-N(\bar{f})}{L^2(0,T)}+\norm{N(\bar{f})}{L^2(0,T)}\nonumber \\
			&\leq \norm{f-\bar{f}}{\X_T}+\norm{N(\bar{f})}{L^2(0,T)} \leq M+\norm{N(\bar{f})}{L^2(0,T)}.
		\end{align}
		If we now combine \eqref{eq:ELipschitz}, \eqref{eq:KappaBounded} and \eqref{eq:N_fBounded}, we obtain from \eqref{eq:ContractionLambdaSplit1}
		\begin{align*}
			&C(f_0)\norm{ \abs{N(f)} \abs{\CalE(f)-\CalE(\tilde{f})}\abs{\vKap_{{f}}}}{L^2(0,T;L^2)} \\
			&\qquad \leq \left(M+\norm{N(\bar{f})}{L^2(0,T)}\right) \blue{C(\bar{f})}\norm{f-\tilde{f}}{\X_T}\leq \frac{q}{4}\norm{f-\tilde{f}}{\X_T}
		\end{align*}
		if we take $\blue{T=T(q,\bar{f}), M=M(q,\bar{f})\in(0,1]}$ small enough.
		
		For the second term in \eqref{eq:ContractionLambda0}, using \eqref{eq:N_f contraction} and \eqref{eq:KappaBounded} we have
		\begin{align*}
			\begin{split}
				C(f_0)\norm{N(f)(t)-N(\tilde{f})(t)}{L^2(0,T)}\norm{\vKap_f(t,x)}{L^{\infty}(0,T;L^2)}&\leq \blue{C(\bar{f})}q_2\norm{f-\tilde{f}}{\X_T} \\
				&\leq \frac{q}{4}\norm{f-\tilde{f}}{\X_T},
			\end{split}
		\end{align*}
		taking $\blue{q_2=q_2(q, \bar{f})\in (0,1)}$ small enough and possibly
		reducing $\blue{T=T(q,\bar{f})}$, $\blue{M=M(q, \bar{f})\in(0,1]}$ if necessary.
		
		For the last term in \eqref{eq:ContractionLambda0}, using \Cref{lem:lowerBoundsEnergy}, \eqref{eq:N_fBounded} and \eqref{eq:KappaLipschitz}, we have
		\begin{align*}
			&\Norm{\frac{N(f)}{\CalE(f)}}{L^{2}(0,T)}\norm{\vKap_f-\vKap_{\tilde{f}}}{L^{\infty}(0,T;L^2)} \\
			&\quad\leq \frac{3}{\CalE(f_0)} \left( M+\norm{N(\bar{f})}{L^2(0,T)}\right) \blue{C(\bar{f})}\norm{f-\tilde{f}}{\X_T}
			\leq \frac{q}{4}\norm{f-\tilde{f}}{\X_T}.
		\end{align*}
		taking $\blue{M=M(q,\bar{f}),T=T(q,\bar{f})\in(0,1]}$ small enough. All in all, we have now shown that taking $\blue{T=T(q,\bar{f}), M=M(q,\bar{f})\in (0,1]}$ small enough, we have
		\begin{align*}
			\norm{\Lambda(f)-\Lambda(\tilde{f})}{L^{2}(0,T;L^2)}=\norm{\lambda(f)\vKap_f-\lambda(\tilde{f})\vKap_{\tilde{f}}}{L^2(0,T;L^2)}\leq \frac{3q}{4}\norm{f-\tilde{f}}{\X_T},
		\end{align*}
		which proves, that $\Lambda\colon\blue{\bar{B}_{T,M}}\to L^2(0,T;L^2)$ is a $\frac{3q}{4}$-contraction. Reducing $\blue{T=T(q,\bar{f})}\in (0,1]$, $\blue{M=M(q,\bar{f}) \in (0,1]}$ if necessary, we may assume that $\CalF$ is a $\frac{q}{4}$-contraction, hence $\mathcal{N}\colon\blue{\bar{B}_{T,M}}\to \Y_T^{1}$ is a $q$-contraction for $\blue{T=T(q, \bar{f}),M=M(q,\bar{f})\in(0,1]}$ small enough.
	\end{proof}

	\section{A gluing lemma for reparametrizations}
	
	In \Cref{thm:LTE}, we used the fact that two smooth reparametrizations can be interpolated by another smooth reparametrization. We state this {gluing} result here in a slightly more general form for possible future reference.
	
	\begin{lem}\label{lem:Repara Interpolation}
		Let $0<t_1<t_2<T$ and $\Phi_1\colon [0,t_2]\times I\to I$, $\Phi_2\colon[t_1, T]\times I\to I$ be smooth families of reparametrizations, such that $\Phi_i(t, \cdot)$ is strictly  increasing for all suitable $t$ and $i=1,2$.
		Then, there exists a smooth family of strictly  increasing reparametrizations $\Psi\colon [0,T]\times I\to I$ satisfying
		\begin{align*}
			\Psi(t,x) &= \Phi_1(t,x), \quad \text{ for all }0\leq t\leq t_1, x\in I\\		
			\Psi(t,x) &= \Phi_2(t,x), \quad \text{ for all }t_2\leq t\leq T, x\in I.
		\end{align*}
	\end{lem} 
	\begin{proof}
		Let $\delta>0$ be sufficiently small and  $\eta\colon [0,T]\to \R, 0\leq \eta\leq 1$ be a smooth cutoff function, satisfying 
		\begin{align*}
			\eta (t) = \left\lbrace \begin{array}{ll}
			1, &\text{ for } 0\leq t\leq t_1+\delta \\
			0, &\text{ for } t\geq t_2-\delta.
			\end{array}\right.
		\end{align*}
		Then it is not difficult to check that the function $\Psi\colon[0,T]\times I\to \R$ given by
		\begin{align*}
			\Psi(t,x) \defeq \left\lbrace\begin{array}{ll}
			\Phi_1(t,x)& \text{ for }0\leq t\leq t_1, x\in I\\
			\Phi_1(t,x) \eta(t) + \Phi_2(t,x)(1-\eta(t)), &\text{ for } t\in [t_1, t_2], x\in I\\
			\Phi_2(t,x) & \text{ for }t_2\leq t\leq T, x\in I
			\end{array}\right.
		\end{align*}
		is smooth and satisfies all the desired properties.
	\end{proof}

\section*{Acknowledgments}
Fabian Rupp has been supported by the Deutsche Forschungsgemeinschaft (DFG,
German Research Foundation)-Projektnummer: 404870139 and by the Austrian Science Fund (FWF), grant numbers \href{https://doi.org/10.55776/P32788}{10.55776/P32788} and \href{https://doi.org/10.55776/ESP557}{10.55776/ESP557}. The authors would like to thank Anna Dall’Acqua, Marius M\"uller, and Rico Zacher for helpful discussions and comments. \blue{ Moreover, the authors are grateful to the referees
	for their valuable feedback on the original manuscript.}

\bibliographystyle{abbrv}
\bibliography{Lib}
\end{document}